\documentclass[12pt]{article}
\usepackage{amssymb,amsthm,amsmath,amsfonts,latexsym,tikz,hyperref,shuffle,ytableau}
\usepackage[hmargin=.8in,vmargin=.6in]{geometry}

\newcommand{\ben}{\begin{enumerate}}
\newcommand{\een}{\end{enumerate}}
\newcommand{\ble}{\begin{lem}}
\newcommand{\ele}{\end{lem}}
\newcommand{\bth}{\begin{thm}}
\renewcommand{\eth}{\end{thm}}
\newcommand{\bpr}{\begin{prop}}
\newcommand{\epr}{\end{prop}}
\newcommand{\bco}{\begin{cor}}
\newcommand{\eco}{\end{cor}}
\newcommand{\bcon}{\begin{conj}}
\newcommand{\econ}{\end{conj}}
\newcommand{\bde}{\begin{defn}}
\newcommand{\ede}{\end{defn}}
\newcommand{\bex}{\begin{exa}}
\newcommand{\eex}{\end{exa}}
\newcommand{\barr}{\begin{array}}
\newcommand{\earr}{\end{array}}
\newcommand{\btab}{\begin{tabular}}
\newcommand{\etab}{\end{tabular}}
\newcommand{\beq}{\begin{equation}}
\newcommand{\eeq}{\end{equation}}
\newcommand{\bea}{\begin{eqnarray*}}
\newcommand{\eea}{\end{eqnarray*}}
\newcommand{\bal}{\begin{align*}}
\newcommand{\bce}{\begin{center}}
\newcommand{\ece}{\end{center}}
\newcommand{\bpi}{\begin{picture}}
\newcommand{\epi}{\end{picture}}
\newcommand{\bpp}{\begin{picture}}
\newcommand{\epp}{\end{picture}}
\newcommand{\bfi}{\begin{figure} \begin{center}}
\newcommand{\efi}{\end{center} \end{figure}}
\newcommand{\bprf}{\begin{proof}}
\newcommand{\eprf}{\end{proof}\medskip}
\newcommand{\capt}{\caption}
\newcommand{\bsl}{\begin{slide}{}}
\newcommand{\esl}{\end{slide}}
\newcommand{\bfr}{\begin{frame}}
\newcommand{\efr}{\end{frame}}

\newcommand{\comp}{\models}

\newcommand{\hqed}{\hfill \qed}

\newcommand{\ol}{\overline}

\newcommand{\hs}[1]{\hspace{#1}}
\newcommand{\hso}[1]{\hspace{-1pt}}
\newcommand{\vs}[1]{\vspace{#1}}
\newcommand{\qmq}[1]{\quad\mbox{#1}\quad}

\newcommand{\emp}{\emptyset}

\newcommand{\sbe}{\subseteq}

\newcommand{\iso}{\cong}

\newcommand{\case}[4]{\left\{\barr{ll}#1&\mbox{#2}\\#3&\mbox{#4}\earr\right.}

\def\<{\langle}
\def\>{\rangle}

\newcommand{\ree}[1]{(\ref{#1})}

\newcommand{\ra}{\rightarrow}

\newcommand{\al}{\alpha}
\newcommand{\be}{\beta}
\newcommand{\ga}{\gamma}
\newcommand{\de}{\delta}
\newcommand{\ep}{\epsilon}

\newcommand{\io}{\iota}

\newcommand{\la}{\lambda}

\newcommand{\si}{\sigma}

\newcommand{\De}{\Delta}

\newcommand{\bx}{{\bf x}}

\newcommand{\bP}{{\bf P}}
\newcommand{\bQ}{{\bf Q}}
\newcommand{\bR}{{\bf R}}

\newcommand{\bbF}{{\mathbb F}}

\newcommand{\bbP}{{\mathbb P}}

\newcommand{\bbZ}{{\mathbb Z}}

\newcommand{\cG}{{\cal G}}

\newcommand{\cS}{{\cal S}}
\newcommand{\cT}{{\cal T}}

\newcommand{\fS}{{\mathfrak S}}

\newcommand{\ab}{\ol{a}}
\newcommand{\Ab}{\ol{A}}

\DeclareMathOperator{\Des}{Des}
\DeclareMathOperator{\Det}{Det}

\DeclareMathOperator{\NSym}{NSym}

\DeclareMathOperator{\sgn}{sgn}

\DeclareMathOperator{\sh}{sh}

\DeclareMathOperator{\st}{st}

\DeclareMathOperator{\Sym}{Sym}

\DeclareMathOperator{\co}{co}
\DeclareMathOperator{\QSym}{QSym}
\DeclareMathOperator{\mQSym}{{\mathfrak m}QSym}
\DeclareMathOperator{\SSym}{\fS Sym}
\DeclareMathOperator{\PSym}{{\mathfrak P} Sym}
\newcommand{\qsh}{\hs{3pt}\ol{\shuffle}\hs{3pt}}
\DeclareMathOperator{\rev}{rev}
\newcommand{\empt}{\rule{0pt}{0pt}}
\newcommand{\shu}{\shuffle}
\newcommand{\Ft}{\tilde{F}}

\newtheorem{thm}{Theorem}[section]
\newtheorem{prop}[thm]{Proposition}
\newtheorem{cor}[thm]{Corollary}
\newtheorem{lem}[thm]{Lemma}
\newtheorem{conj}[thm]{Conjecture}
\newtheorem{exa}[thm]{Example}

\begin{document}
\pagestyle{plain}

\title{Antipodes and involutions
}
\author{Carolina Benedetti\\[-5pt]
\small Department of Mathematics, Michigan State University,\\[-5pt]
\small East Lansing, MI 48824-1027, USA, {\tt caro.benedetti@gmail.com}\\[5pt]
and\\[5pt]
Bruce E. Sagan\\[-5pt]
\small Department of Mathematics, Michigan State University,\\[-5pt]
\small East Lansing, MI 48824-1027, USA, {\tt sagan@math.msu.edu}
}

\date{\today\\[10pt]
	\begin{flushleft}
	                                       \\[5pt]
	\end{flushleft}}

\maketitle

\begin{abstract}

If $H$ is a connected, graded Hopf algebra, then Takeuchi's formula can be used to compute its antipode.  However, there is usually massive cancellation in the result.  We show how sign-reversing involutions can sometimes be used to obtain  cancellation-free formulas.  We apply this idea to nine different examples.  We rederive known formulas for the antipodes in the Hopf algebra of polynomials, the shuffle Hopf algebra, the Hopf algebra of quasisymmetric functions in both the monomial and fundamental bases,  the Hopf algebra of multi-quasisymmetric functions in the fundamental basis, and the incidence Hopf algebra of graphs.  We also find cancellation-free expressions for particular values of the antipode in the immaculate basis for the noncommutative symmetric functions as well as the Malvenuto-Reutenauer and Poirier-Reutenauer Hopf algebras, some of which are the first of their kind.  We include various conjectures and suggestions for future research.

\end{abstract}

%
%

\section{Introduction}

Let $(H,m,u,\De,\ep)$ be a bialgebra over a field $\bbF$.  Call  $H$ {\em graded} if it can be written as $H=\oplus_{n\ge0} H_n$ so that
\ben
\item $H_i H_j \sbe H_{i+j}$ for all $i,j\ge0$,
\item $\De H_n\sbe \oplus_{i+j=n} H_i\otimes H_j$ for all $n\ge0$, and
\item $\ep H_n = 0$ for all $n\ge1$.
\een
If $H_0\iso \bbF$, then we say that $H$ is {\em connected}.  Takeuchi~\cite{tak:fhg}  showed that if a bialgebra is  graded and connected,  then it is a Hopf algebra and gave an explicit formula for its antipode.  To state his result, define a projection map $\pi:H\ra H$ by linearly extending
\beq
\label{proj}
\pi|_{H_n}=\case{0}{if $n=0$,}{I}{if $n\ge1$,}
\eeq
where $0$ and $I$ are the zero and identity maps, respectively.
\bth[\cite{tak:fhg}]
Let $H$ be a  connected graded bialgebra.  Then $H$ is a Hopf algebra with antipode
\beq
\label{tak}
S=\sum_{k\ge0} (-1)^k m^{k-1} \pi^{\otimes k}\De^{k-1},
\eeq
where we let $m^{-1}=u$ and $\De^{-1}=\ep$.\hqed
\eth

Equation~\ref{tak} has the advantage of giving an explicit formula for the antipode.  But it is usually not the most efficient way to calculate $S$ as there is massive cancellation in the alternating sum.  One of the standard combinatorial techniques for eliminating cancellations is the use of sign-reversing involutions.  Let $A$ be a set and $\io:A\ra A$ be an involution on $A$ so that $\io$ is composed of fixed points and two-cycles.  Suppose that $A$ is equipped with a sign function $\sgn:A\ra\{+1,-1\}$.  The involution $\io$ is {\em sign reversing} if, for each two-cycle $(a,b)$, we have 
$\sgn a = -\sgn b$.  It follows that
$$
\sum_{a\in A}\sgn a=\sum_{a\in F} \sgn a,
$$
where $F$ is the set of fixed points of $\io$.  Furthermore, if all elements of $F$ have the same sign, then 
$$
\sum_{a\in A}\sgn a=\pm|F|,
$$
where the bars denote cardinality.

The purpose of the current work is to show how sign-reversing involutions can be used to give cancellation-free formulas for graded connected Hopf algebras.  We give nine different examples of this technique.  The first, in Section~\ref{pha}, is an application to the Hopf algebra of polynomials.  Of course, it is easy to derive the formula for $S$ in this case by other means.  But the ideas of splitting and merging which will appear over and over again can be seen here in their simplest form.  Next, we consider the shuffle Hopf algebra where, again, splitting and merging provide a simple proof.  One also sees why applying $S$ yields the reversed word as it appears naturally as the unique fixed point of our involution.   More complicated applications appear in Sections~\ref{QSymM} and~\ref{QSymF} where we consider the Hopf algebra of quasisymmetric functions in the monomial and fundamental bases.  Motivated by ideas from $K$-theory, Lam and Pylyavskyy~\cite{lp:cha} defined multi-analogues of several Hopf algebras.  Very recently, Patrias~\cite{Pat15:afc} derived cancellation-free expressions for their antipodes and we show how our method can be used to obtain one of them in Section~\ref{mQSym}.  Next we give an involution proof of a formula of Humpert and Martin~\cite{hm:iha} for the antipode in the incidence Hopf algebra on graphs.  Again, the acyclic orientations which are counted by the coefficients appear naturally when finding the fixed points of the involution.  We end with three examples involving the immaculate basis of the Hopf algebra of noncommutative symmetric functions defined by Berg et al.~\cite{bbssz:lsh}, the Malvenuto-Reutenauer Hopf algebra of permutations~\cite{mr:sph} and the Poirier-Reutenauer Hopf algebra of Young tableaux~\cite{pr:aht}.  Some of these expressions are the first cancellation-free ones in the literature. Aguiar and Mahajan~\cite{am:hmha} provided a cancellation-free antipode formula for the Malvenuto-Reutenauer Hopf algebra using Hopf monoids. We recover some of their results using certain involutions, and appealing only to the Hopf algebra structure. We end with a section about future research and open problems, as well as noting other recent work where our technique has been applied.

We should mention that various other researchers have been studying cancellation-free formulae of antipodes.  For example, M\'endez and Liendo~\cite{ml:afn} have constructed a Hopf algebra associated with any symmetric set operad.  They then give a combinatorial formula for its antipode using Schr\"oder trees.  In another direction, Menous and Patras~\cite{MP15:rha} generalize the  forest formula for computing the antipode of the Hopf algebras of Feynman diagrams in perturbative quantum field theory, showing that it can be used in arbitrary right-handed polynomial Hopf algebras.

\section{The polynomial Hopf algebra}
\label{pha}

In this section we will use a sign-reversing involution to derive the well-known formula for the antipode in the polynomial Hopf algebra $\bbF[x]$.  We need  some combinatorial preliminaries.  If $n$ is a nonnegative integer, then let $[n]=\{1,2,\dots,n\}$.  An {\em ordered set partition} of $[n]$ is a sequence of nonempty disjoint subsets $\pi=(B_1,B_2,\dots,B_k)$ such that 
$\uplus_i B_i =[n]$ where $\uplus$ is disjoint union.  We denote this relation by $\pi\models[n]$.  The $B_i$ are called {\em blocks} and since they are sets we are free to always list their elements in a canonical order which will be increasing.  We will also usually leave out the curly brackets and the commas within each block, although we will retain the commas separating the blocks.  So, for example $(13,2)$ has blocks $B_1=\{1,3\}$ and $B_2=\{2\}$; the partition $(2,3,1)$ has blocks $B_1=\{2\}$, $B_2=\{3\}$, $B_3=\{1\}$; and $(123)$ has a single block $\{1,2,3\}$. Finally, it will sometimes be convenient to allow some of the $B_i$ to be empty, in which case we will write $\pi\models_0 [n]$.  

Since the antipode is linear, it suffices to know its action on a basis.  Here we use the standard basis for $\bbF[x]$.

\bth
\label{F[x]}
In $\bbF[x]$ we have
$$
S(x^n)=(-1)^n x^n.
$$
\eth
\bprf
To apply Takeuchi's formula, we first need to describe $\De^{k-1}(x^n)$.
By definition
$$
\De (x) = 1\otimes x + x \otimes 1 = \sum_{(B_1,B_2)\models_0 [1]} x^{|B_1|}\otimes x^{|B_2|}.
$$
It follows from coassociativity that
$$
\De^{k-1} (x) =\sum_{(B_1,\dots,B_k)\models_0 [1]} x^{|B_1|}\otimes\dots\otimes x^{|B_k|},
$$
and since $\De$ is an algebra map
$$
\De^{k-1}(x^n)=(\De^{k-1} x)^n=\sum_{(B_1,\dots,B_k)\models_0 [n]} x^{|B_1|}\otimes\dots\otimes x^{|B_k|}.
$$
Plugging this into equation~\ref{tak} and remembering that $\pi$ kills anything in $H_0$ gives
\beq
\label{S(x^n)}
S(x^n) = \sum_{k\ge0} (-1)^k \sum_{(B_1,\dots,B_k)\models [n]} x^{|B_1|}\dots x^{|B_k|}.
\eeq
Since $|B_1|+\dots+|B_k|=n$ whenever $(B_1,\dots,B_k)\models [n]$, the previous equation simplifies to
$$
S(x^n) =x^n \sum_{k\ge0}\ \sum_{(B_1,\dots,B_k)\models [n]}  (-1)^k.
$$

The last displayed equation shows that we will be done if we can find a sign-reversing involution $\io$ on the set 
$$
A=\bigcup_{k\ge0} \{\pi=(B_1,\dots,B_k)\ :\ \pi\models [n]\},
$$
where
$$
\sgn (B_1,\dots,B_k)=(-1)^k,
$$
and $\io$ has a single fixed point which is 
\beq
\label{phi}
\phi=(n,n-1,\dots,1).
\eeq
This involution will be built out of two other maps.  If $\pi=(B_1,\dots,B_k)$, then the result of {\em merging} blocks $B_i$ and $B_{i+1}$ is the ordered partition $\mu_i(\pi)$ obtained by replacing these two blocks by $B_i\cup B_{i+1}$.  For example, if $\pi=(5,3,249,16,78)$, then $\mu_3(\pi)=(5,3,12469,78)$.  If $|B_i|\ge2$, then the result of {\em splitting} 
$B_i=a_1\dots a_j$  (where, as usual, the elements of the block are listed in increasing order) is the ordered partition $\si_i(\pi)$ obtained by replacing $B_i$ by the ordered pair $a_1,a_2\dots a_j$.  Returning to our example, we have $\si_3(\pi)=(5,3,2,49,16,78)$.

To define $\io(\pi)$ where $\pi=(B_1,\dots,B_k)$ we find the least index $l$, if any, such that either $|B_l|\ge 2$ or 
$B_l=\{a\}$ with $a<\min B_{l+1}$.  If there is such an index, then we let
$$
\io(\pi)=\case{\si_l(\pi)}{if $|B_l|\ge2$,}{\mu_l(\pi)}{else.}
$$
Otherwise $\io(\pi)=\pi$.  Continuing with $\pi=(5,3,249,16,78)$ from the previous paragraph, we can not have $l=1$ or $2$ since $5>3$ and $3>2$.  But $|B_3|\ge2$ which results in $\io(\pi)=\si_3(\pi)=(5,3,2,49,16,78)$.

It is clear from the definition that $\io$ is a sign-reversing map.  To show that $\io^2(\pi)=\pi$, it suffices to consider the case where $\pi$ is not a fixed point.  Given the definition of the index $l$, we must have 
$\pi=(a_1,a_2,\dots,a_{l-1},B_l,\dots,B_k)$ where $a_1>a_2>\dots>a_{l-1}>\min B_l$.
Suppose first that $|B_l|\ge 2$ and let $B_l=a_l a_{l+1}\dots a_m$ as well as $B_l'=a_{l+1}\dots a_m$.  
Then we have 
$$
\io(\pi)=\si_l(\pi)=(a_1,\dots,a_l, B_l',B_{l+1},\dots,B_k).
$$  
Furthermore, 
$$
\io^2(\pi)=\mu_l(a_1,\dots,a_l, B_l',B_{l+1},\dots,B_k))=\pi
$$
because $a_1>\dots>a_l<a_{l+1}=\min B_l'$.  The demonstration that $\io^2(\pi)=\pi$ when $|B_l|=1$ is similar.

There remains to show that the only fixed point of $\io$ is $\phi$ as defined by equation~\ree{phi}.  But if $\pi$ is a fixed point, then the index $l$ does not exist which implies  $\pi=(a_1,a_2,\dots,a_n)$ with $a_1>a_2>\dots>a_n$.  Clearly the only ordered partition of this type is $\pi=\phi$.
\eprf

\section{The shuffle Hopf algebra}

We next use the split-merge technique to derive the antipode of the shuffle Hopf algebra.  Let $A$ be a finite alphabet and consider the Kleene closure $A^*$ of all words $w=a_1\dots a_n$ over $A$.    We let $l(w)=n$ denote the {\em length} of $w$.  The underlying vector space of the {\em shuffle Hopf algebra} is the set of formal sums $\bbF A^*$.  The product is given by shuffling
$$
v\cdot w =\sum_{u\in v\shu w} u,
$$
where $v\shu w$ indicates all $\binom{l(v)+l(w)}{l(v)}$ interleavings of $v$ and $w$.  Note that if different interleavings result in the same final word, then  they are considered distinct and so such words are counted with multiplicity, for example, $ab\cdot a= 2aab+aba$.
The coproduct is
$$
\De w =\sum_{uv=w} u\otimes v,
$$
where $uv$ denotes concatenation, not product, and we permit $u$ or $v$ to be empty.
To state the formula for the antipode we will need, for $w=a_1a_2\dots a_n$, the {\em reversal operator}
$\rev w = a_n a_{n-1}\dots a_1$.

\bth
The antipode in $\bbF A^*$  is given by
$$
S(w) = (-1)^{l(w)}\rev w.
$$
\eth
\bprf
Applying Takeuchi, we see that
\beq
\label{S(w)}
S(w)=\sum_{k\ge1} (-1)^k \sum_{w_1\dots w_k=w} w_1\cdot\ \dots\ \cdot w_k,
\eeq
where none of the $w_i$ are empty.   
Assume that $w=a_1 a_2\dots a_n$ where the $a_i$ are considered distinct variables.  Once we have proved the result for such $w$, the general case will follow by specialization of the $a_i$.  We will use similar reasoning in future proofs without comment.
Because of the distinctness condition,
 $\rev w$ only occurs in the term $w_1\cdot\ \dots\ \cdot w_n$ and does so with sign $(-1)^n$.  So it suffices to give a sign-reversing involution on the rest of the words in the sum.

Let $v\neq\rev w$ be a word resulting as a shuffle in the term $w_1\cdot\ \dots\ \cdot w_k$ of~\ree{S(w)}.  Find the largest index $j\ge0$ such that 
\ben
\item $l(w_1)=\dots=l(w_j)=1$ (which implies $w_i=a_i$ for $i\le j$), and 
\item  $a_j  a_{j-1}\dots  a_1$ is a subword of $v$. 
\een
 Now $a_{j+1}$ is the first letter of $w_{j+1}$.  If $a_j$ is to the left of $a_{j+1}$ in $v$, then $v$ will also be a shuffle of opposite sign in the merged product
$$
a_1\cdot a_2\cdot\ \dots\ \cdot a_{j-1} \cdot a_j w_{j+1} \cdot w_{j+2}\cdot\ \dots\ \cdot w_k.
$$
Our involution will pair these two copies of $v$.  If $a_j$ is to the right of $a_{j+1}$ in $v$, then we must have $l(w_{j+1})>1$ because, if not, then either $j$ would not have been maximum or $v=\rev w$.  Thus $v$ will also be a shuffle  of opposite sign in the split product
$$
a_1\cdot\ a_2\cdot\ \dots\  \cdot a_j \cdot a_{j+1} \cdot w_{j+1}' \cdot w_{j+2}\cdot \dots \cdot w_k,
$$
where $w_{j+1}'$ is $w_{j+1}$ with $a_{j+1}$ removed.  It is easy to see that these two operations are inverses and so we are done.
\eprf

\section{The monomial basis of $\QSym$}
\label{QSymM}

Quasisymmetric functions were introduced by Gessel~\cite{ges:mpi} to study properties of $P$-partitions where $P$ is a poset.  Malvenuto and Reutenauer~\cite{mr:dqf} then showed that the vector space, $\QSym$, of quasisymmertic functions can be given a Hopf algebra structure and that its dual is related to Solomon's decent algebra.
We wish to use involutions to rederive  known formulas for the antipode acting on  two bases for  $\QSym$.  We start with the monomial basis.  
The formula for $S$ in this basis was derived independently by Ehrenborg~\cite{ehr:pha}, and by Malvenuto and Reutenauer~\cite{mr:dqf}.
Here there will turn out to be more than one term in the final sum even though it is cancellation free.  But the split-merge method will show how these summands appear naturally as fixed points of the involution.

Let $\bx=\{x_1,x_2,\dots\}$ be a countably infinite set of variables.  Vector space bases for $\QSym$ are indexed by compositions 
$\al=(\al_1,\al_2,\dots,\al_l)$ which are sequences of positive integers called parts.  The number of parts of $\al$ is called its {\em length} and denoted $l(\al)$.  The {\em monomial quasisymmetric function} corresponding to $\al$ is defined by
$$
M_\al =M_\al(\bx)=\sum_{i_1<i_2<\dots < i_l} x_{i_1}^{\al_1} x_{i_2}^{\al_2}\dots x_{i_l}^{\al_l}.
$$
The $M_\al$ form a basis for $\QSym$.  The product in $\QSym$ is the normal product of power series.  The coproduct is given by
$$
\De M_{\al} = \sum_{\be\ga=\al} M_{\be} \otimes M_{\ga},
$$
where $\be\ga$ is concatenation with the same conventions as in the shuffle algebra.  Applying Takeuchi's formula, we obtain
\beq
\label{S(M)}
S(M_\al)=\sum_{k\ge 1} (-1)^k \sum_{\al_1\dots\al_k=\al} M_{\al_1}\cdot M_{\al_2} \cdot\ \dots\ \cdot M_{\al_k},
\eeq
where all the $\al_i$ are nonempty.
We refer to the terms corresponding to a given index $k$ in the inner sum as the {\em $k$th summand} of $S(M_\al)$.  We will use the notation
\beq
\label{Mprod}
M_{\al_1} \cdot M_{\al_2}\cdot\ \dots\ \cdot M_{\al_k}
=\al_1\cdot\al_2\cdot\ \dots\ \cdot \al_k.
\eeq
Note that, again,  ``$\cdot$" is being used to distinguish multiplication from concatenation.

Suppose that $\al$ has length $l(\al)=l$.  Our strategy will be to cancel all terms in the $k$th summand of~\ree{S(M)} into terms from either the $(k-1)$st or $(k+1)$st summand for $k<l$.  A term from the $l$th summand will either cancel with one from the $(l-1)$st summand or be a fixed point.  We first need to characterize the terms which can occur in the product~\ree{Mprod}.  To do this, we need to recall the notion of a quasishuffle.  Let $A$ be a set of variables and let $v$ be a vector whose components are sums of the variables.  Given $B\sbe A$ then the {\em restriction of $v$ to $B$} is the vector $v|_B$ obtained by setting the variables not in $B$ equal to zero and eliminating any components which are completely zeroed out in this way.  For example, if $v=(c+d+e,b+f,a+c)$ and $B=\{c,d\}$, then $v|_B=(c+d,c)$.   
Now consider compositions $\al$ and $\be$   as two sets of distinct variables. 
In this context, their  {\em quasishuffle} is a vector $v$ containing only these variables such that $v|_\al=\al$ and $v|_\be=\be$.  We let
$$
\al\qsh\be=\{v\ :\ \text{$v$ is a quasishuffle of $\al$ and $\be$}\}.
$$
For example
$$
(a,b)\qsh (c) =\{(a,b,c), (a,b+c), (a,c,b), (a+c,b), (c,a,b)\}.
$$
It is well known that
$$
\al_1\cdot\al_2\cdot\ldots\cdot\al_k =\sum_{v\in \al_1\qsh\al_2\qsh\dots\qsh\al_k} v.
$$

To state the formula for the antipode, we need two more notions.  If $\al=(a_1,a_2,\dots,a_l)$ is a composition, then its {\em reversal} is the composition
$$
\rev\al=(a_l,a_{l-1}\dots,a_1),
$$
just as for words.
We will also use the refinement partial order on compositions.  Define $\be\ge\al$ to mean $\be$ is a coarsening of $\al$, that is, the parts of $\be$ are obtained by adding together adjacent parts of $\al$.
\bth[\cite{ehr:pha},\cite{mr:dqf}]
\label{S(M)thm}
The antipode in the monomial basis of $\QSym$ is given by
$$
S(M_\al)=(-1)^{l(\al)}\sum_{\be\ge\rev(\al)} M_{\be}.
$$
\eth
\bprf
Let $l=l(\al)$.  We first define the action of the splitting operator.  
It will be convenient to define $\si$ on pairs $(\pi, v)$
where  $v$ is a term in the product $\pi=\al_1\cdot\ \dots\ \cdot\al_k$ and $k<l$.  Since $k<l$, there must be an index $j$ with $l(\al_j)\ge2$.  Let $j$ be the smallest such index, so that $l(\al_1)=\dots=l(\al_{j-1})=1$.  The splitting operator is then defined to be $\si(\pi,v)=(\pi',v)$ where
\beq
\label{Msplit}
\pi'=\al_1\cdot\ldots\cdot\al_{j-1}\cdot(d)\cdot\al_j'\cdot\al_{j+1}\cdot\ldots\cdot\al_k,
\eeq
$d$ is the first element of $\al_j$, and $\al_j'$ is $\al_j$ with $d$ removed.    Note that this is well defined since if $v$ is a term in $\pi$, then it must also be a term in $\pi'$ because the variables in $\al_j'$ are a subset of the ones in $\al_j$.
For example, suppose $\pi=(a)\cdot (b)\cdot(c,d,e)\cdot(f,g)$ which contains the term $v=(c+f,d,a+b+e+g)$.  In this case, $j=3$ and $\pi'=(a)\cdot (b)\cdot(c)\cdot(d,e)\cdot(f,g)$.  Note that $\si$ can be applied to any pair $(\pi,v)$ in a product with $k<l$.

The set of pairs to which we can apply the merge map $\mu$ is more restricted.  
For the rest of the proof, we assume that $\al$ has its component variables listed in lexicographic order.
Consider a pair $(\pi,v)$ where $\pi$ has the form~\ree{Mprod}.  In order for $\mu$ to be the inverse of $\si$, we can only merge $\al_j$ and $\al_{j+1}$ if $l(\al_1)=\dots=\l(\al_j)=1$.  Given a quasishuffle $v$ in the product $\pi$, find the largest index $j$, if any, satisfying 
\ben
\item[(i)] $l(\al_1)=\dots=l(\al_j)=1$, and
\item[(ii)] if $B=\al_j\al_{j+1}$, then $v|_B=(d,\dots,e)$ where the elements $d,\dots,e$ are  listed in increasing lexicographic order (as are the elements of $\al$).
\een
Finally, we define $\mu(\pi,v)=(\pi',v)$ where
$$
\pi'=\al_1\cdot\ldots\cdot\al_{j-1}\cdot\al_j\al_{j+1}\cdot\al_{j+2}\cdot\ldots\cdot\al_k.
$$
Note that  condition (ii) and the fact that $v$ is a term in $\pi$ imply  $v$ is also a term in $\pi'$ and so the map is well defined.    To illustrate, suppose that we have $(\pi,v)$ where $\pi=(a)\cdot (b)\cdot(c)\cdot(d,e)\cdot(f,g)$  and  $v=(c+f,d,a+b+e+g)$.  We can not have $j=1$ because then $B=\{a,b\}$ and $v|_B=(a+b)$ with $a$ and $b$ in the same component.  Similarly $j=2$ will not work since then $B=\{b,c\}$ and $v|_B=(c,b)$ which is not in lexicographic order.  But $j=3$ satisfies both conditions, giving $\mu(\pi,v)=(\pi',v)$ where $\pi'=(a)\cdot (b)\cdot(c,d,e)\cdot(f,g)$.  So $\mu$ inverts the action of $\si$ in the previous example.

To define the involution $\io$, take a pair $(\pi,v)$ and define $\io(\pi,v)=\mu(\pi,v)$ if an index satisfying (i) and (ii) above can be found.   If there is no index, then there are two possibilities.  One is that there is an index $j$ such that (i) is true, but $l(\al_{j+1})\ge2$ and $d$, the sole element of $\al_j$ is in the same component as or to the right of the leftmost element of $\al_{j+1}$.  In this case we define $\io(\pi,v)=\si(\pi,v)$.  The minimality of $j$ implies that $\io$ as defined thus far is an involution, and it is clearly sign reversing.

The only other possibility is that there is no index $j$ satisfying (i) and (ii) because $l(\al_i)=1$ for all $i$ and (ii) is never true.  In this case we must have $k=l$, so that $v$ is a quasishuffle from the last summand with sign $(-1)^l$.  We claim that in this case we have that (ii) is not satisfied if and only if  $v\geq\rev(\al)$.  Indeed $v\not\geq\rev(\al)$ is equivalent to the existence of a pair of consecutive letters $B=\{d,e\}$ appearing in lexicographic order in $v$.  Thus $v|_B=(d,e)$ which is the same as saying that (ii) will  be satisfied for some index.  So the fixed points of the involution give us exactly the quasishuffles we need for the $M_\be$ with $\be\ge\rev(\al)$.
\eprf

\section{The fundamental basis for $\QSym$}
\label{QSymF}

The {\em fundamental quasisymmetric function} corresponding to a partition $\al$ can be defined as
$$
F_\al= \sum_{\be\le\al} M_\be.
$$
The formula for $S$ in the fundamental basis first occurs in~\cite{mr:dqf}.  Our proof for the formula for the antipode of $\QSym$ in the fundamental basis will be very similar to the one for the monomial basis.  We will arrange the notation and exposition to emphasize this fact.  

Associated with any composition $\al=(a_1,\dots,a_l)$ is its {\em rim-hook diagram} which has $a_i$ cells in the $i$th row from the bottom and the last cell of row $a_i$ is in the same column as the first cell of row $a_{i+1}$.  We make no distinction between a composition and its diagram.  For example
$$
\al=(3,1,3,2)=\ydiagram{4+2,2+3,2+1,0+3}\ .
$$
A {\em cut-edge} of $\al$ is an edge which is the first vertical edge of $\al_1$, the last vertical edge of $\al_l$, or an edge bounding two cells of $\al$.  Separating $\al$ into pieces along a cut-edge results in two diagrams $\be$ to the southwest and $\ga$ to the northeast.  In this case we write $\al=\be|\ga$.  The coproduct applied to a fundamental quasisymmetric function is then
$$
\De F_\al=\sum_{\be|\ga=\al} F_\be\otimes F_\ga,
$$
To illustrate, if $\al=(3,1)$ then the various pairs for $\be$ and $\ga$ are, as the cut-edge travels from southwest to northeast,
$$
\emp \hs{15pt} \ydiagram{2+1,3}\ 
, \hs{20pt}
\begin{ytableau}
\none & \none & \none & \empt \\
\empt & \none & \empt & \empt
\end{ytableau}\
, \hs{20pt}
\begin{ytableau}
\none & \none & \none & \empt \\
\empt & \empt & \none & \empt
\end{ytableau}\
, \hs{20pt}
\begin{ytableau}
\none & \none & \empt \\
\none & \none & \none \\
\empt & \empt & \empt
\end{ytableau}\
, \hs{20pt}
\ydiagram{2+1,3} \hs{15pt} \emp\ .
$$
Thus
$$
\De(F_{(3,1)})= 1\otimes F_{(3,1)} + F_{(1)}\otimes F_{(2,1)} + F_{(2)}\otimes F_{(1,1)} + F_{(3)}\otimes F_{(1)}
+F_{(3,1)}\otimes 1.
$$

To describe the product of fundamental quasisymmetric functions we  associate compositions with permutations.  Write $\al=(\al_1,\dots,\al_l)\models n$ 
or $|\al|=n$
if $\sum_i\al_i=n$.  There is a canonical bijection between $\al\models n$ and subsets of $[n-1]$ given by sending $\al$ to 
\beq
\label{D(al)}
D(\al)=\{\al_1,\al_1+\al_2,\dots,\al_1+\al_2+\dots+\al_{l-1}\}.
\eeq
Given a  sequence of integers $w=c_1 c_2\dots c_n$ we denote its {\em descent set} by
\beq
\label{Des}
\Des w = \{i\ :\ c_i>c_{i+1}\}.
\eeq
So every such  $w$ has an associated set $\Des w\sbe [n-1]$ which corresponds to a composition $\al$.  In this case we define $F_w=F_\al$ and say that $w$ {\em models} $\al$.  Note that one can tell by context whether the subscript is a word $w$ or a composition $\al$ since the latter will  have parentheses and commas while the former will not.
Now suppose $\al\models m$ and $\be\models n$.  Let $w_\al$  and $w_\be$ be disjoint (as sets) and model $\al$ and $\be$, repectively.  In this situation, the multiplication of fundamental quasisymmetric functions is given by
\beq
\label{Fprod}
F_\al F_\be = \sum_{w\in w_\al\shu w_\be} F_w,
\eeq
where the sum is over all ordinary shuffles of $w_\al$ and $w_\be$.

Given $\al$ and a set $C$ of positive integers with $|\al|=C$, there is a canonical way to construct a $w$ modeling $\al$ with entries in $C$.  Fill the cells of the diagram of $\al$ bijectively with the elements of $C$ so that rows and columns increase to form a tableau $T$.  So, in particular, if $C=[n]$, then $T$ is a standard Young tableau of shape $\al$.  The {\em row word} of $T$, $w_T$, is constructed by concatenating the rows of $T$ starting with the bottom row and moving up.  Continuing the example started at the beginning of this section, we could take the tableau in Figure~\ref{568413927}
in which case $w_T=568413927$.  It is easy to see that the row and column restrictions on $T$ imply that $w_T$ models $\al$.

\begin{figure}
$$
T=\begin{ytableau}
\none & \none & \none & \none & 2 & 7\\
\none & \none &    1     &      3    & 9\\
\none & \none &    4\\
   5     &     6    &    8\\
\end{ytableau}
$$
\capt{An example for computing the row word}
\label{568413927}
\end{figure}

We are now ready to put everything together and apply Takeuchi's formula.  Fix once and for all a standard Young tableau $T$ of shape $\al$.  Then
\beq
\label{S(F)eq}
S(F_\al)=\sum_{k\ge 0} (-1)^k \sum_{\al_1|\dots|\al_k=\al} F_{w_{T_1}}\dots F_{w_{T_k}},
\eeq
where the $\al_i$ are all nonempty and $T_i$ is the subtableau cut out from $T$ by $\al_i$.  To illustrate, suppose $\al=(2,1)$ and 
$$
T=
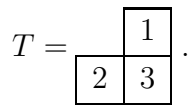
\begin{ytableau}
\none & 1\\
   2     & 3
\end{ytableau}\ .
$$
Thus the terms in $S(F_{(2,1)})$ correspond to the decompositions
$$
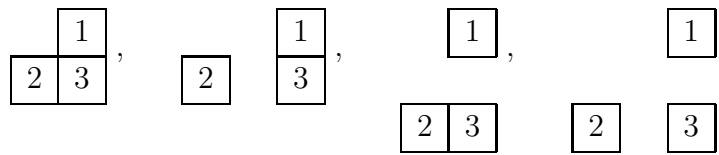
\begin{ytableau}
\none & 1\\
   2     & 3
\end{ytableau}\
, \hs{20pt}
\begin{ytableau}
\none & \none & 1\\
   2     & \none & 3
\end{ytableau}\
, \hs{20pt}
\begin{ytableau}
\none & 1\\
\none & \none\\
   2     & 3
\end{ytableau}\
, \hs{20pt}
\begin{ytableau}
\none & \none & 1\\
\none & \none & \none\\
   2     & \none & 3
\end{ytableau}
$$
of $T$ so that
$$
S(F_{(2,1)})= -F_{231} + F_2 F_{31} + F_{23} F_1 - F_2 F_3 F_1.
$$

We will adopt the terminology and notation of the previous section, the only differences being that we will use products of words to stand for products of fundamental quasisymmetric functions and that such products will be ordinary shuffles because of~\ree{Fprod}.  We will use the notation $\al^t$ to stand for the {\em transpose} or {\em conjugate} of $\al$, that is, the diagram obtained by reflecting $\al$ in the main diagonal.
\bth[\cite{mr:dqf}]
\label{S(F)thm}
The antipode in the fundamental basis of $\QSym$ is given by
$$
S(F_\al)=(-1)^{|\al|} F_{\al^t}.
$$
\eth
\bprf
Let $n=|\al|$.  The definition of the splitting map is  the same as in the proof of Theorem~\ref{S(M)thm} except that one is dealing with products $\pi$ of words and shuffles $v$ which are terms in $\pi$.  Also, in this case, we can apply the splitting map to any product with $k$ factors where $k<n$ since $n$ is the maximum number of nonempty subcompositions into which $\al$ can be decomposed by cuts.  To illustrate, suppose $\al=(2,3,1)$ and we fix the tableau in Figure~\ref{562341}.  
Now the product $\pi=5 \cdot 6 2 3 4\cdot 1$ will contain the shuffle 
$v= 6 2 1 3 4 5$.  Since $j=2$ is the smallest (in fact, only) index of a word in the product of length larger than one, we will have $\si(\pi,v)=(\pi',v)$ where
$\pi'=5 \cdot 6 \cdot 2 3  4\cdot 1$.

\begin{figure}
$$
T=\begin{ytableau}
\none & \none & \none & 1\\
\none &    2     &     3     & 4\\
    5    &    6
\end{ytableau}
$$
\capt{An example for $F_\al$
\label{562341}}
\end{figure}

The merge map is again very similar to the one used for the monomial quasisymmetric functions.  We consider a shuffle $v$ in the product $\pi=w_1\cdot\ldots\cdot w_k$ and find the smallest index $j$, if any, such that
\ben
\item[(i)] $|w_1|=\dots=|w_j|=1$, and
\item[(ii)] if $B=w_j w_{j+1}$, then $v|_B$ is the row word of a rim-hook subtableau of $T$.
\een
We then let $\mu(\pi,v)=(\pi',v)$ where
$$
\pi'=w_1\cdot\ \dots\ \cdot w_{j-1} \cdot w_j w_{j+1} \cdot w_{j+2} \cdot\ \dots\ \cdot w_k
$$
which is well defined, as before, because of condition (ii).
Taking as an example $\pi=5 \cdot 6 \cdot 2 3  4\cdot 1$ and $v= 6 2 1 3 4 5$ we can not have $j=1$ since then $B=\{5,6\}$ and $v|_B=65$ which is not the row word of the first two squares of $T$.
 But when $j=2$ we have $x|_B=6234$ which is the row word of the four middle squares of $T$.  Thus $\mu(\pi,v)=(\pi',x)$ where 
$\pi'=5 \cdot 6 2 3 4\cdot 1$, undoing the previous example's application of $\si$.

One now defines the involution exactly as was done for the $M$-basis, applying $\si$ if possible and otherwise appying $\mu$ if possible.  So the only thing new is to determine the fixed point(s) $(\pi,v)$.  As before, they will all be in the last summand of equation~\ree{S(F)eq} which corresponds to $k=n$.  So the factors in $\pi$ are the individual elements of $T$ and there is only one choice for $\pi$.  Furthermore if $\mu$ can not be applied to $(\pi,v)$, then one must have every pair of adjacent elements in $w_T$ being in the reverse order in $v$.  But this can only happen if $v$ is $w_T$ read backwards, which is precisely the row word of $\al^t$.  
\eprf

\section{The fundamental basis of $\mQSym$}
\label{mQSym}

Motivated by work of Buch~\cite{buc:lrr} on set-valued tableaux, Lam and Pylyavskyy~\cite{lp:cha} defined six new Hopf algebras. These can be thought of as $K$-theoretic analogues of the symme\-tric function, quasisymmetric function, noncommutative symmetric function, and Malvenuto-Reutenauer Hopf algebras (the first and last both having two analogues).  They appealed to Takeuchi's Theorem to conclude the existence of antipodes.  Recently, Patrias~\cite{Pat15:afc} has given explicit formulas for these maps.  We wish to show how one of these expressions can be derived using splitting and merging.

We will describe multi-$\QSym$, denoted $\mQSym$, the $K$-theoretic analogue of $\QSym$.  Because of the use of set-valued maps, we will need to permit arbitrary $\bbZ$-linear combinations of basis elements and these elements will not be homogeneous of a certain degree.    Since elements of $\mQSym$ are not of bounded degree, it is not graded.  However, we can still apply Takeuchi's formula because its proof also works more generally for any Hopf algebra  where the projection map~\ree{proj} is locally nilpotent.

Let $\tilde{\bbP}$ be the family of all finite, nonempty sets of positive integers.  If $S,T\in\tilde{\bbP}$, then we write $S<T$ (respectively, $S\le T$) if $\max S < \min T$ (respectively, $\max S\le \min T$).  Let $w=c_1 c_2\dots c_n$ be a permutation of $[n]$.  A {\em $w$-set-valued partition} is a map $f:[n]\ra\tilde{\bbP}$ satisfying
$$
\case{f(i)\le f(i+1)}{if $c_i<c_{i+1}$,}{f(i)< f(i+1)}{if $c_i>c_{i+1}$}
$$
for $i\in[n-1]$.
This is a special case of a more general definition for $P$-set-valued partitions, $P$ a poset, which we will not need.  To illustrate, if $w=231$, then a $w$-set-valued partition would satisfy $f(1)\le f(2)< f(3)$.  For example one could have $f(1)=\{5,7\}$, $f(2)=\{7,8,10\}$ and $f(3)=\{11\}$.  

Associate with any $S\in\tilde{\bbP}$ the monomial $\bx_S=\prod_{s\in S} x_s$, and with any $w$-set-valued partition $f$ the monomial $x_f=\prod_{i\in[n]} x_{f(i)}$.  The {\em fundamental multi-quasisymmetric function} associated with a composition $\al$ is
$$
\Ft_\al=\Ft_w =\sum_f x_f,
$$ 
where $w$ models $\al$, and the sum is over all $w$-set-valued partitions $f$.  Continuing the example above,  $\al=(2,1)$ is modeled by $w=231$ and the given partition would contribute $x_5 x_7^2 x_8 x_{10} x_{11}$ to $\Ft_{(2,1)}$.  We note that the sum of the terms of least degree in $\Ft_\al$ is exactly $F_\al$.  Finally we let $\mQSym$ be the span of the $\Ft_\al$.

We need some combinatorial constructions to describe the bialgebra structure of $\mQSym$.
Many of the ideas which came into play in proving the antipode formula in the fundamental basis for $\QSym$ will also be used here.
In addition to being able to separate  a diagram at a cut-edge, we will need to be able to separate it at a cell.  So if $c$ is a cell of the diagram of $\al$ which we will call the {\em cut-cell}, then we write $\al=\be\bullet\ga$ where $\be$ is the composition whose diagram is all cells southwest of and including $c$ and $\ga$ is the composition to the northeast and including $c$.  
Equivalently, $\al$ is formed by identifying the last square of $\be$ with the first square of $\ga$.
Note that both $\be$ and $\ga$ include $c$ so that
$|\be|+|\ga|=|\al|+1$.  For example, if $\al=(3,1)$, then here are the various pairs $\be$ and $\ga$ as the cell $c$ moves from southeast to northwest
$$
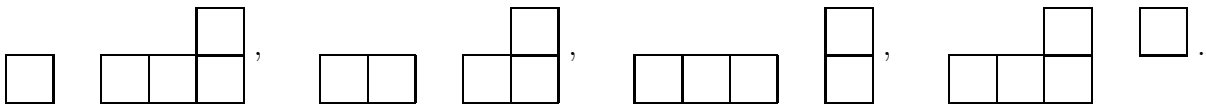
\begin{ytableau}
\none & \none & \none & \none & \empt \\
\empt & \none & \empt & \empt & \empt
\end{ytableau}\
, \hs{20pt}
\begin{ytableau}
\none & \none & \none & \none & \empt \\
\empt & \empt & \none & \empt  & \empt
\end{ytableau}\
, \hs{20pt}
\begin{ytableau}
\none & \none  & \none & \none & \empt \\
\empt & \empt  & \empt & \none & \empt
\end{ytableau}\
, \hs{20pt}
\begin{ytableau}
\none & \none  & \empt & \none & \empt \\
\empt & \empt  & \empt & \none & \none
\end{ytableau}\
.
$$
The coproduct  of $\mQSym$ can now be written
$$
\De(\Ft_\al) =\sum_{\be,\ga} \Ft_\be \otimes \Ft_\ga,
$$
where the sum is over all $\be,\ga$ such that $\al=\be|\ga$ or $\al=\be\bullet\ga$.  Continuing our example
\begin{align*}
\De(\Ft_{(3,1)})& =1\otimes \Ft_{(3,1)} + \Ft_{(1)}\otimes \Ft_{(3,1)}
+\Ft_{(1)}\otimes \Ft_{(2,1)} + \Ft_{(2)}\otimes \Ft_{(2,1)}
+ \Ft_{(2)}\otimes \Ft_{(1,1)}\\ 
&\qquad\qquad\qquad+  \Ft_{(3)}\otimes \Ft_{(1,1)}
+ \Ft_{(3)}\otimes \Ft_{(1)} + \Ft_{(3,1)}\otimes \Ft_{(1)}
+\Ft_{(3,1)}\otimes 1
\end{align*}
as the position of the cut travels over alternating edges and cells from southwest to northeast.

We can now apply Takeuchi's formula to get
\beq
\label{S(Ft)}
S(\Ft_\al)=\sum_{k\ge0} (-1)^k \sum_{\al_1,\dots,\al_k} \Ft_{\al_1}\dots\Ft_{\al_k}
\eeq
with the sum being over all $\al_1,\dots,\al_k$ such that
\beq
\label{decomp}
\al=\al_1\circ_1\al_2\circ_2\dots\circ_{k-1}\al_k,
\eeq
where either $\circ_i=|$ or $\circ_i=\bullet$ for all $i$.  Note that two expressions with the same $\al_i$ but different $\circ_i$ both contribute separately to the sum.  For example, if $\al=(3)$, then $\al=(1)\bullet (2)|(1)$ and $\al=(1)|(2)\bullet (1)$ are different terms.  We can write any expression of the form~\ree{decomp} by grouping together all the compositions between any two occurrences of an edge cut.  Specifcially, we will write
$$
\al_1\circ_1\al_2\circ_2\dots\circ_{k-1}\al_k = \be_1|\be_2| \dots | \be_m,
$$
where $\be_1=\al_1\bullet\dots\bullet\al_a$, $\be_2=\al_{a+1}\bullet\dots\bullet\al_b$, and so forth.  We will call the $\be_i$ {\em components} and the $\al_j$ {\em subcomponents} of this expression.  To illustrate, if $\al=(3,1,1)$, then
$$
(3,1,1) = (1)\bullet(2)\bullet(1) | (1) | (1)\bullet(1,1)
$$
has three components, namely $\be_1 =  (1)\bullet(2)\bullet(1), \be_2=(1), \be_3= (1)\bullet(1,1)$.

It will be convenient to have a geometric way to visualize the components and subcomponents of an expression.  For the components, we will use the same convention as for $\QSym$, where the diagram of $\al$ is split along the cut-edges.  For the subcomponents, each cut-cell will be split in two by a vertical edge.  These edges will be decorated with a bullet to distinguish them from the original edges of $\al$.  The diagram for the example at the end of the last paragraph is shown in Figure~\ref{DecFig}.
Note that the components of the expression are the connected components of the diagram while the subcomponents can be obtained by cutting each component along the bullet edges.

\begin{figure}
\bce
\begin{ytableau}
\none & \none & \none & \none & \none & \none & \empt\\
\none & \none & \none & \none & \none & \hs{19pt}\bullet & \empt\\
\none\\
\hs{19pt}\bullet & \empt &\hs{19pt}\bullet & \empt & \none & \empt
\end{ytableau}
\hs{60pt}
\begin{ytableau}
\none & \none & \none & \none & \none & \none & 1\\
\none & \none & \none & \none & \none & \hs{11pt} 2 \hs{2pt}\bullet & 3\\
\none\\
\hs{11pt} 4 \hs{2pt}\bullet & 5 & \hs{11pt} 6 \hs{2pt}\bullet & 7 & \none & 8 
\end{ytableau}
\ece
\capt{The diagram of $(1)\bullet(2)\bullet(1) | (1) | (1)\bullet(1,1)$ and it superstandard labeling
\label{DecFig}
}
\end{figure}

To complete our exposition, we need to describe how to take products of the $\Ft_\al$.  Given a word $w=c_1\dots c_n$, then a {\em multiword} on $w$ is a word of the form $\tilde{w}=c_1^{m_1} \dots c_n^{m_n}$ where $c_i^{m_i}$ indicates that $c_i$ is to be repeated $m_i>0$ times for all $i$.  Exponents of one can be ommited.  For example, if $w=231$, then we could take $\tilde{w}=2^4 3 1^2=2222311$.  The {\em multishuffles} of words $v$ and $w$, denoted $v\tilde{\shu} w$, is the set of all words $x=d_1\dots d_r$ such that
\begin{enumerate}
\item $x$ is a shuffle of some $\tilde{v}$ and $\tilde{w}$, and
\item $d_i\neq d_{i+1}$ for all $i\in[r-1]$.
\end{enumerate}
Note that this set is infinite.
By way of illustration if $v=21$ and $w=3$, then
$$
21\tilde{\shu} 3 = \{213,\ 231,\ 321,\ 2321,\ 2131,\ 2313,\ 3213,\ 3231,\ \dots\},
$$
where we have listed all the multishuffles of length $3$ or $4$.  Now given compostions $\al$ and $\be$ we take words $w_\al$ and $w_\be$ modeling them, respectively, on disjoint alphabets.  In this case
\beq
\label{Ftprod}
\Ft_\al \Ft_\be = \sum_{w\in w_\al \tilde{\shu} w_\be} \Ft_w.
\eeq
Continuing our example, $21$ models $\al=(1,1)$ and $3$ models $\be=(1)$ so that
$$
\Ft_{(1,1)} \Ft_{(1)}=\Ft_{213}+\Ft_{231}+\Ft_{321}+\Ft_{2321}+\Ft_{2131}+\Ft_{2313}+\Ft_{3213}+\Ft_{3231}+\dots\ .
$$

To combine~\ree{S(Ft)} and~\ree{Ftprod}, consider a decomposition of the form~\ree{decomp}.  Label the corresponding diagram to form a tableau in a superstandard way with the numbers $1,\dots,\sum_i|\al_i|$ from left to right in each row starting with the top row and working down.  Then, for each $i$, let $w_i$ be the row word of the subtableau of $T$ corresponding to $\al_i$.  Our decomposition will be represented by 
\beq
\label{wdecomp}
\pi=w_1\circ_1 \dots \circ_{k-1} w_k
\eeq
and the terms in the corresponding product of multi-quasisymmetric functions will be indexed by the multishuffles in $w_1\tilde{\shu}\dots\tilde{\shu}w_k$.  Figure~\ref{DecFig} illustrates these ideas, showing that  $(1)\bullet(2)\bullet(1) | (1) | (1)\bullet(1,1)$ when written in words becomes $4 \bullet 56 \bullet 7 | 8 |  2 \bullet 31$ which contains all terms of $S(F_{(3,1)})$ corresponding to the multishuffle 
$4\tilde{\shu}56 \tilde{\shu} 7\tilde{\shu} 8\tilde{\shu}  2 \tilde{\shu} 31$.

To describe the coefficients of the cancellation free formula for the antipode we need the notion of collapsing a diagram.  This operation is called merging in~\cite{Pat15:afc}, but that would conflict with our use of the term in this work.    We say that $\be$ is a {\em collapse} of $\al$ if one can succesively collapse together  boxes of $\al$ which share an edge to form $\be$. We let $c_{\al,\be}$ be the number of ways to collapse $\al$ to $\be$.  In counting collapses only the sets of boxes collapsed matters, not the order of the collapsing.  Figure~\ref{col} shows that if $\al=(3,1,1)$ and $\be=(2,1)$, then $c_{\al,\be}=4$.  The labeling of the boxes is merely to show which sets were collapsed.

\begin{figure}
\bce
$\al =$
\begin{ytableau}
\none & \none & 1\\
\none & \none & 2\\
3 & 4 & 5
\end{ytableau}
\hs{40pt}
$\be:$
\begin{ytableau}
\none & 1\\
34 & 25
\end{ytableau}
\hs{20pt}
\begin{ytableau}
\none & 12\\
34 & 5
\end{ytableau}
\hs{20pt}
\begin{ytableau}
\none & 12\\
3 & 45
\end{ytableau}
\hs{20pt}
\begin{ytableau}
\none & 1\\
3 & 245
\end{ytableau}
\ece
\capt{Collapsing $\al=(3,1,1)$ onto $\be=(2,1)$ 
\label{col}}
\end{figure}

\bth[\cite{Pat15:afc}]
The antipode in the  fundamental basis of $\mQSym$ is given by
$$
S(\Ft_\al)=\sum_\be (-1)^{|\be|} c_{\be,\al^t} \Ft_\be,
$$
where the sum is over all compositions $\be$.
\eth
\bprf
The proof parallels that of Theorem~\ref{S(F)thm}.  All diagrams are labeled in  a superstandard way and their reading words used in the corresponding multishuffles.  We will define the involution on pairs $(\pi,v)$ where $\pi$ is a decomposition of the form~\ree{decomp} and $v$ is a term in the multishuffle corresponding to $\pi$.  We denote the image of the pair under the involution as $(\pi',v')$.  As usual, the split operation is easiest to describe.  For all $i$, let  $v_i$ be the subword of $v$ which is a multiword on $w_i$ where $w_i$ is the reading word of $\al_i$ in the superstandard tableau for $\pi$.  Find the smallest  index $j$, if any, such that
$|v_j|\ge 2$ and suppose $v_j=ab\dots c$.  There are now two cases.
\ben
\item If $a\neq b$, then let 
	\ben
	\item $\pi'=\pi$ with $\al_j$ replaced by $(1)|\al_j'$ where $\al_j'$ is $\al_j$ with its first square removed, and
	\item $v'=v$.
	\een
\item If $a = b$, then let 
	\ben
	\item $\pi'=\pi$ with $\al_j$ replaced by $(1)\bullet\al_j$, and
	\item $v'=v$ with one added to all elements greater than or equal to  $b$ except $a$. 
	\een  
\een

\begin{figure}
$$
\begin{ytableau}
\none &  1  & 2\\
\none\\
\hs{11pt} 3 \hs{2pt}\bullet & 4 
\end{ytableau}\ ,
\quad
1342
\qquad
\barr{c}
\\
{\si \atop\longrightarrow}\\[10pt]
{\longleftarrow \atop \mu}
\earr
\qquad
\begin{ytableau}
\none & 1 & \none & 2\\
\none\\
\hs{11pt} 3 \hs{2pt}\bullet   &  4 
\end{ytableau}\ ,
\quad
1342
$$
\vs{10pt}
$$
\begin{ytableau}
\none & 1\\
\none\\
2  & 3
\end{ytableau}\ ,
\quad 
2123
\qquad
\barr{c}
\\
{\si \atop\longrightarrow}\\[10pt]
{\longleftarrow \atop \mu}
\earr
\qquad
\begin{ytableau}
\none & \none  & 1\\
\none \\
 \hs{11pt} 2 \hs{2pt}\bullet & 3  &  4 
\end{ytableau}\ ,
\quad
2134
$$
\capt{The splitting and merging operations
\label{(2,1)}}
\end{figure}

To illustrate these cases, let $\al=(1,2)$.  Suppose  $\pi=(1)\bullet (1)|(2)$ which in terms of words is $3\bullet 4|1 2$ as can be seen in the top line of Figure~\ref{(2,1)}.  Let $v=1342$, a multishuffle in $\pi$ with corresponding subwords $v_1=3$, $v_2=4$, and $v_3=12$.  The only $v_j$ with at least two elements is $12$ and for this subcomponent $1\neq 2$. 
So using the first case above, we have $\si(\pi,v)=(\pi',v')$ where $\pi'=(1)\bullet(1)|(1)|(1)$ and $v'=1342$.
On the other hand, suppose we consider $\al=(2,1)$ and  $\pi=(2)|(1)$ or in term of  words $23|1$ as in the bottom line of Figure~\ref{(2,1)}.  Then $\pi$ contains the multishuffle $v=2123$ with subwords $v_1=223$ and $v_2=1$.  Now $v_1$ is the only subword with at least two elements and it begins with $22$.  So we are in the second case of the defiition of $\si$.  Thus $\pi'=(1)\bullet(2)|(1) $.  Also $v'$ is obtained from $v$ by increasing the second $2$ and all larger numbers by one to obtain $v'=2134$.  Note this convention is precisely what is needed to make $v'$ a shuffle in $\pi'$ and so this case is well defined.  It is even easier to see that the first case is as welll.

To describe the merge map, consider $(\pi,v)$ and find the smallest index $j$, if any, such that $|v_i|=1$ for $i\le j$ and the concatenation $v_j v_{j+1}$ is a multisubword of the concatenation $w_j w_{j+1}$. Furthermore, if $\circ_j=\bullet$ we also insist that $b$ and $c$  are not consecutive in $v$ where $v_j=b$ and $c$ is the first element of $v_{j+1}$.  Again, there are two cases to define $\mu(\pi,v)=(\pi',v')$.
\ben
\item If $\circ_j=|$, then let 
	\ben
	\item $\pi'=\pi$ with $\al_{j}|\al_{j+1}$ replaced by $\al_j'$ which is formed by edge-pasting $\al_{j}$ and $\al_{j+1}$ back together again, and
	\item $v'=v$.
	\een
\item If $\circ_j=\bullet$, then let 
	\ben
	\item $\pi'=\pi$ with $\al_{j}\bullet \al_{j+1}$ replaced by $\al_j'$which is formed by cell-pasting $\al_{j}$ and $\al_{j+1}$ back together again, and
	\item $v'=v$ with one subtracted from all element larger than or equal to  $c$.
	\een  
\een

By way of example, consider $\pi=(1)\bullet(1)|(1)|(1)$ and $v=1342$ which is the upper right pair in Figure~\ref{(2,1)}.  In terms of words $\pi=3\bullet4|1|2$ and $v$ has corresponding subwords $v_1=3$, $v_2=4$, $v_3=1$, $v_4=2$.  First consider $j=1$.  Then $v_1v_2=34$ is a subword of $v$.  But $\circ_1=\bullet$ and $3,4$ are adjacent in $v$ so they can not be merged.  Next we try $j=2$ giving $v_2v_3=41$.  This is not a subword of $v$, so we test $j=3$. Finally $v_3v_4=12$ is a subword of $v$ and $\circ_3=|$ so there is no further restriction.  It follows that we can apply the first case of the merge definition and return to the original pair which started the split example.  Now let us look at $\pi=(1)\bullet(2)|(1) $ and $v=2134$ as in the lower right of  Figure~\ref{(2,1)}.  Translating to words gives  $\pi=2\bullet 34| 1$ and $v_1=2$, $v_2=34$, $v_3=1$.  Taking $j=1$ we see that $v_1v_2=234$ is a subword of $v$.  Also $\circ_1=\bullet$, but $2,3$ are not adjacent in $v$ so that we can merge as in case 2.  This undoes the effect of $\si$ in the second split example.  It should not be hard for the reader to prove that $\si$ is well defined.

We now define the involution $\io$ exactly as in the case of $\QSym$: applying $\si$ if an appropriate index $j$ can be found and otherwise applying $\mu$.  It is straightforward to verify that the minimality of $j$ makes these two operations inverses.  In fact case 1 for $\si$ is the inverse of case 1 for $\mu$ and similarly for the 2nd cases.  And, as usual, the fact that $\io$ is sign-reversing is trivial.

There remains to determine the fixed points of $\io$.  For this we will need a notion which is also useful in the theory of permutation patterns.  Call a permutation $v$ of the interval $[a,e]$ {\em colayered} if it is of the form 
$$
v=  d\ (d+1)\dots e\ c\ (c+1)\  \dots (d-1)\ \dots a\ (a+1) \dots (b-1)
$$
for certain $a,b,c\dots,e$.  The reason for this terminology is because the complement of $v$ is layered in the usual sense.
The {\em  layer lengths} of $v$ are the cardinalities of the maximal increasing intervals.   For example $v=678945123$ is colayered with layer lengths $4,2,3$.  There is a natural bijection between colayered permutations of $[d]$ and compositions of $\al=(\al_1,\dots,\al_k)$ of $d$: just assign to each colayered permutation the composition of its layer lengths.  For the inverse, take the diagraram of $\al$ and let $v$ be the row word of its superstandard filling.  So as $v$ varies over all colayered permutations of $[d]$ the associated  $\al$ will run over all compositions of $d$.

\begin{figure}
$$
\begin{ytableau}
\none & \none & \none & \none & \hs{11pt} 1 \hs{2pt}\bullet & 2\\
\none\\
 \hs{11pt} 3 \hs{2pt}\bullet & \hs{11pt} 4 \hs{2pt}\bullet & 5 & \none & 6
\end{ytableau}\ ,
\quad
v=216534
\qquad
\longleftrightarrow
\qquad
\begin{ytableau}
\none &534\\
21 & 6
\end{ytableau}
\quad
\text{ collapsed from }
\quad
\begin{ytableau}
\none & 3 & 4\\
\none & 5\\
1 & 6\\
2
\end{ytableau}
$$
\capt{A fixed point for $\al=(2,1)$}
\label{fp}
\end{figure}

No let $(\pi,v)$ be a fixed point of $\io$ where $\pi=\al_1\circ_1\dots\circ_{k-1} \al_k$.  
Such a fixed point for $\al=(2,1)$ is shown in Figure~\ref{fp} to help clarify the following argument.
Since we can not apply $\si$, we must have $|v_i|=1$ for all $i$.  This has several consequences.  First of all $v$ can have no repeated elements.  Also,  $|\al_i|=1$ for all $i$.  Thus every cut-edge of $\al$ (except the first and last) has been split to form $\pi$ and  there is a natural bijection between the cells of $\al$ and the components of $\pi$.  Furthermore,  every edge internal to a component is dotted.  Let the reading words of the components of $\pi$ be $w_1,\dots,w_l$ in order from southwest to northeast.  In the example, $w_1=345$, $w_2=6$, and $w_3=12$.   

I claim that the  possible second coordinates in our fixed point are are exactly those of the form $v=w_l' w_{l-1}' \dots w_1'$ where $w_i'$ is a colayered permutation of the elements of $w_i$.  And once this claim is established, we will be done.  Indeed,  recall  the bijection between colayered permutations and compositions as well as the bijection between components of $\pi$ and the cells of $\al$. This combined with the order reversal in going from the $w_i$ to the $w_i'$ results in a bijection between the fixed points and row words of compositions $\be$ which collapse to $\al^t$ where the elements of $v_i'$ collapse to form the $i$th box (counting from southeast to northwest) of $\al^t$.  See Figure~\ref{fp} for an illustration.  

We will prove the claim in the case that the number of components is $l=2$ since the general case is similar inducting on $l$.  So let $w_1=a_1\dots a_p$ and $w_2=b_1\dots b_q$. 
Since the components are single cells, the subcomponents of a given component form a single row.
It follows that the reading words $w_1,w_2$ are increasing sequences of consecutive integers.
Since we can not apply $\mu$ we have, by case 1, that $b_1$ must be to the left of $a_p$ in $v$.  And by case 2, the only way that $a_{i+1}$ can be to the right of $a_i$ is if they are adjacent, with a similar statement holding for the $b_j$.  It follows by induction on $i$ and $j$ that all the $b_j$ must come before all the $a_i$.  A similar induction shows that the $a_i$ and $b_j$ must be arranged in colayered permutations.  This completes the proof.
\eprf

\section{The incidence Hopf algebra on graphs}
\label{hag}

Now we turn our attention to the incidence Hopf algebra on graphs $\mathcal G$.  We begin by introducing  some notation. Let $G=(V,E)$ be a simple graph with vertex set $V=[n]$ and edge set $E$.  We denote by $[G]$ the isomorphism class of $G$.
The Hopf algebra $\mathcal G$ is the free $\bbF$-module on the isomorphism classes $[G]$. It has been studied by Schmitt~\cite{sch:hac}, and Humpert and Martin~\cite{hm:iha} among others. Ardila and Aguiar
(private communication) have recently described the antipode of $\mathcal G$ using the more general setting of Hopf monoids. 

The product and coproduct maps in $\mathcal G$ are  as follows
\begin{align}\label{prod}
[G\cdot H] &= [G\uplus H]\\
\label{coprod}\Delta\left([G]\right) &= \sum_{(V_1,V_2)\models [n]}[G_{V_1}]\otimes [G_{V_2}],
\end{align}
where $G_{V_1}$ is the subgraph of $G$ induced by the vertex set $V_1$ and similarly for $G_{V_2}$.
We will henceforth drop the square brackets denoting equivalence class to simplify notation.  No confusion should be caused by this convention.

To state the formula of Humpert and Martin for the antipode, we recall that a {\em flat} of $G$ is a spanning subgraph $F$ such that every component of $F$ is induced.  These are just the flats of the cycle matroid of $G$.  We will denote the number of components of $F$ by $c(F)$.   Given any collection $F$ of edges of $G$, the graph obtained from $G$ by contracting the edges in $F$ will be denoted $G/F$.   Finally, we let $a(G/F)$ be the number of acyclic orientations $O$ of $G/F$.  That is, $O$ is a digraph with underlying graph $G/F$ such that $O$ contains no directed cycles.
\bth[\cite{hm:iha}]
\label{hm}
The antipode in $\mathcal G$ is given by
$$
S(G) = \sum_{F\text{ flat of  }G}(-1)^{c(F)}a(G/F)F.
$$
\eth
\begin{proof}
We will  give  a combinatorial  proof  of this theorem  using the ideas from Section \ref{pha}. 
By virtue of~\ree{prod}, \ree{coprod}, and Takeuchi's formula, the same reasoning that lead to equation~\ree{S(x^n)} applied to $\cG$ gives
\beq
\label{S(G)}
S(G) = \sum_{k\geq 0}(-1)^k\sum_{(V_1,\dots,V_k)\models [n]}G_{V_1}\uplus\cdots\uplus G_{V_k}.
\eeq
Since each $G_{V_i}$ is induced, each graph $F=G_{V_1}\uplus\cdots\uplus G_{V_k}$  in the inner sum is a flat of $G$.
In order to show that the coefficient of $F$ in $S(G)$ is $(-1)^{c(F)} a(G/F)$ we will construct, for each flat $F$, a sign-reversing involution $\iota$ on the  set
$$
A_F= \bigcup_{k\ge0} \{\pi=(V_1,\dots,V_k)\ :\ \pi\models [n]\text{ and }G_{V_1}\uplus\cdots\uplus G_{V_k} = F\}
$$
with sign function
$$
\sgn \pi=(-1)^k
$$
for $\pi=(V_1,\dots,V_k)$.  This involution will have fixed points which are in bijection with acyclic orientations of $G/F$ and which all have sign $(-1)^{c(F)}$, so we will be done.

We will first give the proof in the case that $F$ is the flat with no edges and then indicate how this demonstration can be modified for a general flat.   Note that if $\pi=(V_1,\dots,V_k)\in A_F$ where $F$ is the empty flat, then each $V_i$ is a set of independent vertices. 
Let $O_{\pi}$ be the orientation in $G=G/F$ defined by orienting its edges as
$$
u\rightarrow v\text{ whenever }u\in V_i,\, v\in V_j\, \text{ and }i<j.
$$
Note that $O_{\pi}$ is acyclic because each $V_i$ is independent and the usual ordering of the integers is transitive.  

In general, there will be many different $\pi$ giving rise to the same orientation.  So given an orientation $O$, we let
$$
\Pi_O = \{\pi\ :\ O_\pi=O\}.
$$
Among all the partitions in $\Pi_O$, we distinguish a canonical one $\phi_O=(b_1^O,b_2^O,\dots,b_n^O)$ defined as follows.  Let $b_1^O$ be the largest source in $O$.  Next consider the digraph $O-b_1^O$ and let $b_2^O$ be its largest source.  Continue in this way, removing sources and selecting the largest source in the remaining graph until we are left with a single vertex $b_n^O$.  An example follows this proof.  Our strategy will be, for each acyclic orientation $O$, to construct a sign-reversing involution $\io_O:\Pi_O\ra\Pi_O$  which has $\phi_O$ as its unique fixed point.  Note that $\sgn\phi_O=(-1)^n=(-1)^{c(F)}$ where $F$ is the empty flat.  This will complete the proof in the case under consideration.

To define $\io=\io_O$, suppose $\pi=(V_1,\dots,V_k)\in\Pi_O$ with $\pi\neq\phi_O$.  Then there must be a smallest vertex $i$ where $\pi$ and $\phi_O$ disagree, that is, $V_1=b_1^O,\dots, V_{i-1}=b_{i-1}^O$, but $V_i\neq b_i^O$.  For simplicity, let $b=b_i^O$.  Let $V_j$, $j\ge i$, be the block containing $b$.
If $|V_j|\ge2$, then we {\em split} $V_j$ by replacing it with the ordered pair of blocks $V_j-b,b$ to obtain a partition $\si_j(\pi)$.  If $|V_j|=1$, then $i\neq j$ since otherwise $\pi$ and $\phi_O$ would not differ at index $i$. Consider the suborientation $O'$ obtained by restricting $O$ to $V_i\uplus\dots\uplus V_k$.  So $V_{j-1}$ is in $O'$   which forces $V_{j-1}\uplus V_j$ to be independent: if there were an edge in the graph underlying $O'$ from a vertex of $V_{j-1}$ to $V_j=b$, then it would be oriented into $b$, contradicting the fact that $b$ is a source in $O'$.  So we can {\em merge} $V_j$ by replacing it with $V_{j-1}\uplus V_j$ to form $\mu_j(\pi)$.  Finally define
$$
\io(\pi)=\case{\si_j(\pi)}{if $|V_j|\ge2$,}{\mu_j(\pi)}{if $|V_j|=1$.}
$$

As usual, it is clear that $\io$ is sign-reversing.  To check that it is an involution, we first show that the indices $i$ for $\pi$ and for $\io(\pi)$ are the same.  If a block is split, then $b$ ends up one block to the right of its original position so that the $i$th block still differs from $\phi_O$.  If two blocks are merged, then the only way to change $V_i$ is by merging in with $V_{i+1}$.  But afterwards the $i$th block has size at least two and so again differs from that block in $\phi_O$.

Since $i$ does not change in passing from $\pi$ to $\io(\pi)$, the orientation $O'$ must be the same for both.  Thus $b$ is also invariant under this map since it is the largest source in $O'$.  Now the fact that $\io$ is an involution follows in much that same way as in the proof for $\bbF[x]$.
 This finishes the demonstration for the empty flat.

We now deal with the general case.  As already noted, every term in the sum~\ree{S(G)} is a flat of $G$. Consider a flat $F$ and a term $G_{V_1}\uplus\dots\uplus G_{V_k}=F$.  Then $G/F$ is obtained by contracting every component of $F$ to a point.  Each such component lies in one of the $G_{V_i}$.  So the partition $(V_1,\dots,V_k)$ of $V(G)$ induces a partition $(W_1,\dots,W_l)$ of $V(G/F)$ where each of the $W_i$ is independent in $G/F$.  Clearly this process is reversible with each partition of $V(G/F)$ into independent sets giving rise to a partition of $V(G)$ which induces the flat $F$.  Now we can apply the same process as with the empty flat to the partitions of $V(G/F)$.  This completes the proof.
\end{proof}

As an example of the involution in Theorem~\ref{hm}, consider the acyclic orientation $O$ depicted in the first graph of Figure~\ref{hm:fig}.  To compute the fixed point $\phi_O$ we first remove the largest source, which is vertex $5$, to be the first component of $\phi_O$.  The largest source in what remains is $8$ since the arc from $5$ to $8$ has been removed, and this becomes the second component of $\phi_O$.  Continuing in this way, we obtain  
$$
\phi_O=(5,8,7,3,6,4,2,1)
$$
as displayed in the second drawing of Figure~\ref{hm:fig}.

\bfi
\begin{tikzpicture}
\draw(0,0) node{};
\draw(0,1.3) node{};
\end{tikzpicture}

\begin{tikzpicture}
[inner sep=2, vertex/.style={circle,fill}]
\draw(-1.5,1.5) node{$O=$};
\node at (1,0) [vertex] (6) {};
\draw(1,-.5) node{$6$};
\node at (2,0) [vertex] (5) {};
\draw(2,-.5) node{$5$};
\node at (0,1) [vertex] (7) {};
\draw(-.5,1) node{$7$};
\node at (3,1) [vertex] (4) {};
\draw(3.5,1) node{$4$};
\node at (0,2) [vertex] (8) {};
\draw(-.5,2) node{$8$};
\node at (3,2) [vertex] (3) {};
\draw(3.5,2) node{$3$};
\node at (1,3) [vertex] (1) {};
\draw(1,3.5) node{$1$};
\node at (2,3) [vertex] (2) {};
\draw(2,3.5) node{$2$};
\draw[->,thick] (5)--(8);
\draw[->,thick] (8)--(7);
\draw[->,thick] (5) --(6);
\draw[->,thick] (8)--(1);
\draw[->,thick] (3)--(6);
\draw[->,thick] (3)--(4);
\draw[->,thick] (2)--(1);
\end{tikzpicture}

\vs{-30pt}

\begin{tikzpicture}
[inner sep=2, vertex/.style={circle,fill}]
\draw(-1.8,0) node{$\phi_O=(5,8,7,3,6,4,2,1)=$};
\node at (1,0) [vertex] (5) {};
\draw(1,-.5) node{$5$};
\node at (2,0) [vertex] (8) {};
\draw(2,-.5) node{$8$};
\node at (3,0) [vertex] (7) {};
\draw(3,-.5) node{$7$};
\node at (4,0) [vertex] (3) {};
\draw(4,-.5) node{$3$};
\node at (5,0) [vertex] (6) {};
\draw(5,-.5) node{$6$};
\node at (6,0) [vertex] (4) {};
\draw(6,-.5) node{$4$};
\node at (7,0) [vertex] (2) {};
\draw(7,-.5) node{$2$};
\node at (8,0) [vertex] (1) {};
\draw(8,-.5) node{$1$};
\draw[->,thick] (5)--(8);
\draw[->,thick] (8)--(7);
\draw[->,thick] (5) .. controls +(right:3pt) and +(up:40pt) .. (6);
\draw[->,thick] (8) .. controls +(right:3pt) and +(up:70pt) .. (1);
\draw[->,thick] (3)--(6);
\draw[->,thick] (3) .. controls +(right:1pt) and +(up:20pt)  .. (4);
\draw[->,thick] (2)--(1);
\end{tikzpicture}

\vs{-10pt}

\begin{tikzpicture}
[inner sep=2, vertex/.style={circle,fill}]
\draw(-1.8,0) node{$\pi=(5,3,4,26,8,7,1)=$};
\node at (1,0) [vertex] (5) {};
\draw(1,-.5) node{$5$};
\node at (2,0) [vertex] (3) {};
\draw(2,-.5) node{$3$};
\node at (3,0) [vertex] (4) {};
\draw(3,-.5) node{$4$};
\node at (4,1) [vertex] (6) {};
\draw(4.4,1) node{$6$};
\node at (4,-1) [vertex] (2) {};
\draw(3.6,-1) node{$2$};
\node at (5,0) [vertex] (8) {};
\draw(5,-.5) node{$8$};
\node at (6,0) [vertex] (7) {};
\draw(6,-.5) node{$7$};
\node at (7,0) [vertex] (1) {};
\draw(7,-.5) node{$1$};
\draw[->,thick] (5) .. controls +(right:1pt) and +(up:30pt)  .. (8);
\draw[->,thick] (8)--(7);
\draw[->,thick] (5) .. controls +(right:3pt) and +(up:40pt) .. (6);
\draw[->,thick] (8) .. controls +(right:1pt) and +(up:20pt) .. (1);
\draw[->,thick] (3) .. controls +(right:1pt) and +(down:20pt)  ..(6);
\draw[->,thick] (3)-- (4);
\draw[->,thick] (2) .. controls +(right:1pt) and +(down:20pt)  ..(1);
\end{tikzpicture}

\vs{-10pt}

\begin{tikzpicture}
[inner sep=2, vertex/.style={circle,fill}]
\draw(-1.8,0) node{$\io(\pi)=(5,3,4,268,7,1)=$};
\node at (1,0) [vertex] (5) {};
\draw(1,-.5) node{$5$};
\node at (2,0) [vertex] (3) {};
\draw(2,-.5) node{$3$};
\node at (3,0) [vertex] (4) {};
\draw(3,-.5) node{$4$};
\node at (4,1) [vertex] (6) {};
\draw(4.4,1) node{$6$};
\node at (4,-1) [vertex] (2) {};
\draw(3.6,-1) node{$2$};
\node at (4,0) [vertex] (8) {};
\draw(4,-.5) node{$8$};
\node at (5,0) [vertex] (7) {};
\draw(5,-.5) node{$7$};
\node at (6,0) [vertex] (1) {};
\draw(6,-.5) node{$1$};
\draw[->,thick] (5) .. controls +(right:1pt) and +(up:30pt)  .. (8);
\draw[->,thick] (8)--(7);
\draw[->,thick] (5) .. controls +(right:3pt) and +(up:40pt) .. (6);
\draw[->,thick] (8) .. controls +(right:1pt) and +(up:20pt) .. (1);
\draw[->,thick] (3) .. controls +(right:1pt) and +(down:20pt)  ..(6);
\draw[->,thick] (3)-- (4);
\draw[->,thick] (2) .. controls +(right:1pt) and +(down:20pt)  ..(1);
\end{tikzpicture}

\capt{Example for Theorem~\ref{hm}
\label{hm:fig}
}
\efi

Now suppose we are given $\pi=(5,3,4,26,8,7,1)$ as in the third illustration of Figure~\ref{hm:fig}.  It is easy to verify that $O_\pi=O$.  Comparing $\pi$ and $\phi_O$, we see that they first differ in the the second block so $i=2$.   The largest source  of $O'=O-5$  is $b=8$ which is in the singleton block $V_5$.  Thus
$$
\io(\pi)=\mu_5(\pi)=(5,3,4,268,7,1)
$$
as drawn at the end of Figure~\ref{hm:fig}.

Finally, consider $\pi'=\io(\pi)=(5,3,4,268,7,1)$.  The reader should have no trouble verifying that for $\pi'$ we still have $i=2$ and $b=8$ which is now in block $4$.  Since this block has other elements in it as well
$$
\io^2(\pi)=\io(\pi')=\si_4(\pi')= (5,3,4,26,8,7,1)=\pi
$$
as desired.

Viewed as maps on ordered partitions, the involutions for $\bbF[x]$ and for $\cG$ are not the same.  However, if we take $G$ to be the graph with $V=[n]$ and no edges, then there is only one flat, namely $F=G$, and only one acyclic orientation $O$.  So the involution in this case has a unique fixed point which is the same as the one in the proof of Theorem~\ref{F[x]}.  In fact,  we could have used this map to prove the polynomial result and emphasize even more the similarly of the demonstrations.  We chose the earlier involution because of its simplicity.

\section{The immaculate basis for $\NSym$}
\label{ibn}

The Hopf algebra of noncommutative symmetric functions, $\NSym$, was introduced by Gelfand et al.~\cite{gkllrt:nsf}, and its immaculate basis, $\cS_\al$, was defined in the paper of Berg et al.~\cite{bbssz:lsh}.
There is no known cancellation-free formula for the antipode acting on $\cS_\al$ for general $\al$.  So in this section we derive  expressions in the special cases where $\al$  is of hook shape or has at most two parts.   In the case of hooks, the proof follows easily from the know expression for the antipode acting on the ribbon Schur basis.  For shapes of at most two parts we give a new, cancellation-free expression proved by applying an involution.

The noncommutative symmetric functions are freely generated as an algebra by the noncommutative symbols $H_n$ for $n$ a positive integer.  It is also convenient to let $H_0=1$ and $H_n=0$ for $n<0$.  Similarly, the ordinary symmetric functions, $\Sym$, are generated by the complete homogeneous symmetric functions $h_1,h_2,\dots$ which do commute.  There is also the forgetful function $\NSym\ra\Sym$ defined by algebraically extending the map $H_n\mapsto h_n$ for all $n$.

Bases for $\Sym$ are indexed by integer partitions $\la=(\la_1,\la_2,\dots,\la_k)$, which are weakly decreasing sequences of positive integers.  The all-important Schur function basis can be defined by the $k\times k$ Jacobi-Trudi determinant
$$
s_\la=\det(h_{\la_i+j-i}).
$$
Bases for $\NSym$ are indexed by compositions $\al=(\al_1,\al_2,\dots,\al_k)$, which are arbitrary sequences of positive integers.
To get a  basis for $\NSym$ corresponding to the Schur functions, we define the noncommutative determinant of a $k\times k$ matrix $A=(a_{i,j})$ to be
$$
\Det A = \sum_{\si} a_{1,\si(1)} a_{2,\si(2)}\dots a_{k,\si(k)},
$$
where the sum is over all permutations $\si$ of $[k]$.  The {\em immaculate basis} for $\NSym$ is defined to be 
$$
\cS_\al=\Det(H_{\al_i+j-i}).
$$
If $\al$ is a partition, then we clearly have $\cS_\al\mapsto s_\al$ under the forgetful map.  To simplify the notation for products in this determinant we use, for any sequence of integers $\al=(\al_1,\al_2,\dots,\al_k)$, the shorthand $H_\al=H_{\al_1} H_{\al_2}\dots H_{\al_k}$.

It is well known that $\Sym$ is actually a Hopf algebra where
$$
\De h_n = \sum_{i=0}^n h_i\otimes h_{n-i}.
$$
The antipode has a particularly nice action on the Schur basis, namely
\beq
\label{S(s)}
S(s_\la)=(-1)^{|\la|} s_{\la^t},
\eeq
where $\la^t$ is the conjugate of $\la$.  We use the notation $|\la|$ and $l(\la)$ as we have done for compositions. 

We also have that $\NSym$ is a Hopf algebra with 
\beq
\label{DeHn}
\De  H_n = \sum_{i=0}^n H_i\otimes H_{n-i}.
\eeq
But it appears as if the antipode is much harder to compute in the immaculate basis.  So we will only derive formulas for it  when $\al$ is a hook or when $\al$ has (at most) two rows.

We will first derive a formula for $S(\cS_\al)$ when $\al$ is a hook.  To do this, we recall another important basis of $\NSym$.  The {\em ribbon Schur function} corresponding to $\al$ is
\beq
\label{R_al}
R_\al=\sum_{\be\ge\al} (-1)^{l(\al)-l(\be)} H_\be.
\eeq
The antipode in the ribbon basis has a simple formula which can be found in the book of Grinberg and Reiner~\cite[Theorem 5.42]{GR14:hac}.
\bth
\label{S(R)}
For any composition $\al$,
$$
S(R_\al)=(-1)^{|\al|} R_{\al^t}
$$
\eth

A {\em hook} is a composition of the form $\al=(n,1^k)$.  In this case, the determinant for $\cS_\al$ can be expanded around its first column since the second entry in that column is $H_0=1$ which commutes with the other $H_i$.  This results in the recursion 
$$
\cS_{n,1^k} =H_n \cS_{1^k}  - \cS_{n+1,1^{k-1}}.
$$
On the other hand, partitioning the terms in the sum~\ree{R_al} for $R_{n,1^k}$ into those with $\be_1=n$ and those with $\be_1>n$ shows that this ribbon Schur function satisfies the same recursion.  So the next result follows easily using Theorem~\ref{S(R)} and  induction on $k$.
We note that Grinberg~\cite{Gri15:dpa} has rederived this formula using his work on double posets, a concept introduced by Malvenuto and Reutenauer in~\cite{mr:sph}.

\bth
\label{hook}
For all $n\ge1$ and $k\ge0$,
$$
S(\cS_{n,1^k})=(-1)^{n+k} \cS_{k+1,1^{n-1}}.
$$
\eth

For the two-row case, we will express the antipode in terms of certain sets of tableaux.  The {\em shape} of a composition 
$\al=(\al_1,\dots,\al_l)$ is an array of $l$  rows of left-justified boxes with $\al_i$ boxes in row $i$.  
We will use English notation where the first row is at the top as well as  matrix coordinates for the cells.  We also do not distinguish between a composition and its shape.  So, for example,
$$
\scalebox{.5}{
\begin{tikzpicture}
\draw (-2.2,2) node{\scalebox{2}{(3,1,2,2) = }};
\draw (0,0) grid (2,1);
\draw (0,1) grid (2,2);
\draw (0,2) grid (1,3);
\draw (0,3) grid (3,4);
\end{tikzpicture}.
}
$$
A {\em dual immaculate tableau} of shape $\al$ is a placement $T$ of positive integers in the cells of $\al$ such that the rows strictly increase and the first column weakly increases.  The reason for using ``dual" is because the strong and weak inequalities are interchanged from those for an immaculate tableau as defined in the paper of Berg et al.~\cite{bbssz:lsh}. 
We write $\sh T=\al$.   One dual immaculate tableau of shape $(3,1,2,2)$ is
$$
T=\barr{lll}
1&3&4\\
1\\
2&6\\
2&4
\earr.
$$
We let $T_c=T_{i,j}$ be the element of $T$ in cell $c=(i,j)$.  The {\em content} of $T$ is the composition $\co(T)=(m_1,m_2,\dots)$ where $m_i$ is the multiplicity of $i$ in $T$.  In our example tableau $\co(T)=(2,2,1,2,0,1)$.

Suppose $\cT$ is a set of tableaux.  A {\em set of frozen  cells} for $\cT$ is a set  of cells such that, for each such cell $c$, the element $T_c$ is the same for all $T\in\cT$.  This includes the case when $T_c$ is empty for all $T\in\cT$.  We will denote a frozen cell by giving its element a star.  In the case the cell is to be empty, we use the symbol $\emp^*$.  To illustrate, here is a set of tableaux indicating one of its sets of frozen cells
$$
\cT=
\left\{
\barr{lll}
1^*&3^*&4\\
1^*&4\\
2^*&6^*&\emp^*\\
2^*\\
\emp^*
\earr,
\qquad
\barr{lll}
1^*&3^*&4\\
1^*\\
2^*&6^*&\emp^*\\
2^*&4\\
\emp^*
\earr,
\qquad
\barr{lll}
1^*&3^*\\
1^*&4\\
2^*&6^*&\emp^*\\
2^*&4\\
\emp^*
\earr
\right\}.
$$
In all cases of interest to us, the set of frozen cells will have the shape of a composition, that is, frozen cells in a row are left-justified and the set of frozen cells in the first column is connected.  So  we will call the elements in these cells a {\em frozen tableaux} $T^*$.  Note that $T^*$ includes the cells which are forced  to be empty and we also require that all such cells are either at the right end of a row or at the bottom of the first column.  In our example, the shape of $T^*$ is $(2,1,3,1,1)$.   A dual immaculate tableau which includes empty cells in this way will be called an {\em extended tableau}.

Now given an extended tableau $T^*$ and a content vector $v$, define the set $\cT(T^*,v)$ to be the set of all dual immaculate tableaux such that
\ben
\item[(a)] $T^*$ is a frozen tableau for $\cT(T^*,v)$,
\item[(b)] $\co(T)=v$ for every $T\in\cT(T^*,v)$, and
\item[(c)] $\cT(T^*,v)$ contains every dual immaculate tableau satisfying (a) and (b).
\een
Note that our example $\cT$ is such a set.  It will turn out that the sets $\cT(T^*,v)$ which we will need always also have the property
\ben
\item[(d)] if $v=(v_1,\ldots,v_m)$, then $\co T^* = (v_1,\dots,v_{m-1},w)$ for some $w\le v_m$.
\een
So henceforth we also assume that $\cT(T^*,v)$ also satisfies (d).
It turns out that these are exactly the sets we need to describe the antipode.  

Given $\al=(\al_1,\dots,\al_l)\comp n$ we will have to compute expressions of the form
\begin{align*}
S(H_{\al_1}  H_{\al_2}\dots H_{\al_l})&=S(H_{\al_l})\dots S(H_{\al_2}) S(H_{\al_1})\\
&=(-1)^n \cS_{1^{\al_l}} \dots  \cS_{1^{\al_2}}  \cS_{1^{\al_1}},
\end{align*}
where the second equality comes from the case $k=0$ of Theorem~\ref{hook}.
In~\cite[Theorem 7.3]{bbssz:lsh} a rule is given for expanding $\cS_\al \cS_\la$ in the immaculate basis whenever $\al$ is a composition and $\la$ is a partition.  Applying this rule repeatedly to the last equation easily gives
\beq
\label{S(H)}
S(H_\al)=(-1)^{|\al|}\sum_{T} \cS_{\sh(T)},
\eeq
where the sum is over all dual immaculate tableaux $T$ with $\co(T)=(\al_l,\dots,\al_2,\al_1)$.
We are now in a position to prove our result for two-rowed tableaux.
An example illustrating the following theorem follows its proof.  One can also obtain a formula for $S(\cS_{m,n})$ by changing to the ribbon basis, using Theorem~\ref{S(R)}, and changing basis back, but this expression in not cancellation free.

\bth\label{2row}
Given  $m,n\ge 1$, let $\cT_1=\cT(T_1^*,(n,m))$ and $\cT_2=\cT(T_2^*,(n-1,m+1))$ where
$$
T_1^*=
\barr{ll}
1^*\\
\vdots\\
1^* &2^*
\earr,
\qquad
T_2^*=
\barr{l}
1^*\\
\vdots\\
1^* \\
\emp^*
\earr.
$$
Then
$$
S(\cS_{m,n})=(-1)^{m+n}\sum_{T\in\cT_1} \cS_{\sh T} +(-1)^{m+n+1} \sum_{T\in\cT_2} \cS_{\sh T}.
$$
\eth
\bprf
Since $\cS_{m,n}=H_{m,n}-H_{m+1,n-1}$ we see, using equation~\ree{S(H)}, that
$$
S(\cS_{m,n})=(-1)^{m+n}\sum_T \cS_{\sh T} + (-1)^{m+n+1}\sum_T \cS_{\sh T},
$$
where the first sum is over all possible dual immaculate tableaux of content $(n,m)$ and the second over such tableaux of content $(n-1,m+1)$.   Let $\cT$ be the signed set which is the union of these two sets of tableaux, with signs being assigned so as to give the sums above.  Now it suffices to find a sign-reversing involution $\io:\cT\ra\cT$ whose fixed points are exactly the tableaux in $\cT_1\cup\cT_2$.

Consider $T\in\cT$.  If $T$ contains $n$ ones, then change the lowest one (which must be in the first column since $T$ is dual immaculate) to a two as long as the resulting tableau is still dual immaculate.  If $T$ contains $n-1$ ones, then change the highest two in the first column of $T$ (if one exists) to a one.  Note that if this is possible, then the resulting tableau must be dual immaculate.  If neither of these options is possible, then $T$ is a fixed point.  

It is clear from the definitions that this is an involution and reverses sign.  To find the fixed points, note that the only tableaux with $n$ ones fixed by $\io$ are those where the lowest one also has a two in its row.  This gives precisely the tableaux in $\cT_1$.  Similarly, the tableaux with $n-1$ ones fixed by $\io$ are exactly those with no two in the first column which correspond to the tableaux in $\cT_2$.  This finishes the proof.
\eprf

By way of illustration, to calculate $S(\cS_{2,4})$ we compute
$$
\cT_1
=\left\{
\barr{ll}
1^*&2\\
1^*\\
1^*\\
1^*&2^*\\
\rule{0pt}{5pt}
\earr,
\qquad
\barr{ll}
1^*\\
1^*&2\\
1^*\\
1^*&2^*\\
\rule{0pt}{5pt}
\earr,
\qquad
\barr{ll}
1^*\\
1^*\\
1^*&2\\
1^*&2^*\\
\rule{0pt}{5pt}
\earr,
\qquad
\barr{ll}
1^*\\
1^*\\
1^*\\
1^*&2^*\\
2
\earr
\right\}
$$
and
$$
\cT_2=\left\{
\barr{ll}
1^*&2\\
1^*&2\\
1^*&2\\
\emp^*
\earr
\right\}
$$
so
$$
S(\cS_{2,4})=\cS_{2,1,1,2}+\cS_{1,2,1,2}+\cS_{1,1,2,2}+\cS_{1,1,1,2,1}-\cS_{2,2,2}.
$$

\section{The Malvenuto-Reutenauer Hopf algebra}

Aguiar and Mahajan~\cite{am:hmha} derived an antipode formula for the Malvenuto-Reutenauer Hopf algebra of permutations using an antipode formula in a certain Hopf monoid.   While their formula is cancellation-free, one needs the monoid structure for its construction.  We will derive certain cancellation-free formuas which can be derived without appealing to the monoid.  In particular, we will find such expressions for permutations whose image under the Robinson-Schensted map is a column superstandard Young tableau of hook shape.  These tableaux will appear again in the next section when we consider the Poirier-Reutenauer Hopf algebra.

Let $\fS_n$ be the symmetric group on $[n]$ and $\fS=\cup_{n\ge0} \fS_n$.  
The Malvenuto-Reutenauer Hopf algebra $\SSym$ has basis $\fS$.  To describe the product, if $\si=a_1 a_2\dots a_n \in\fS_n$ and $m$ is a positive integer, then let $\si+m$ denote the sequence obtained by increasing every element of $\si$ by $m$.  For example, if $\si=231$, then $\si+4=675$.  Now if $\pi\in\fS_m$ and $\si\in\fS_n$, then we define
$$
\pi\cdot\si =\sum_{\tau\in \pi\shu(\si+m)} \tau
$$
To illustrate 
$$
12\cdot 21 = 1243+1423+ 1432+ 4123+4132+4312.
$$
For the coproduct, we need the notion of standardization.  If $\si=a_1a_2\dots a_n$ is any sequence of distinct positive integers, then its {\em standardization} is the permutation $\st(\si)$ obtained by replacing the smallest $a_i$ by one, the next smallest by two, and so on.  By way of illustration $\st(9587) = 4132$.  For $\pi\in\fS_n$ we let
$$
\De(\pi)=\sum_{\si\tau=\pi} \st(\si)\otimes \st(\tau),
$$
where $\si\tau$ represents concatenation of sequences, empty sequences allowed.  As an example
$$
\De(3142) = \ep\otimes 3142+1\otimes 132+21\otimes 21 + 213\otimes 1 + 3142\otimes \ep,
$$
where $\ep$ is the empty permutation.

Each term in the Takeuchi expansion of $S(\si)$  is the sum of the elements of a set of shuffles $\si_1\shu\dots\shu\si_k$ where $\si=\st(\si_1)\dots \st(\si_k)$.  It will be convenient in what follows to identify the shuffle set with the sum of its elements.

The next result permits us to derive information about two  antipode expansions at once.
We write $[\pi]f$ for the coefficient of $\pi$ in any formal sum of permutations $f$.  

\bth
\label{inverse}
If $\pi,\si\in\fS_n$, then
$$
[\pi] S(\si)=[\si^{-1}] S(\pi^{-1}).
$$
\eth
\bprf
There is a bijection between the shuffle sets in $S(\si)$ and compositions $\al$ where the shuffle set $\si_1\shu\dots\shu\si_k$ corresponds to the composition $\al=(\al_1,\dots,\al_k)$  with $|\si_i|=\al_i$ for all $i$.  The sign of the shuffle set  is $(-1)^k$ and $\pi$ occurs at most once in each shuffle set.  So to prove the theorem, it suffices to show that $\pi$ occurs in the shuffle set of $S(\si)$ corresponing to $\al$ if and only if $\si^{-1}$ appears in the shuffle set  of $S(\pi^{-1})$ corresponding to $\al$.  By symmetry, it suffices to show the forward implication.

Suppose that $\pi$ occurs in the shuffle set  $\si_1\shu\dots\shu\si_k$.  Then for all $i$ we must have that $\si_i$ is a subsequence of $\pi$.  We will show that in this case the shuffle set $\pi'_1\shu\dots\shu\pi'_k$ in $S(\pi^{-1})$ corresponding to the same composition must contain a copy of $\si^{-1}$ in that $\pi'_i$ is a subsequence of $\si^{-1}$ for all $i$.
We will do the case $i=1$ as the others are similar.  Suppose $\pi=a_1\dots a_n$, $\si=b_1\dots b_n$ and $\si_1=\st(b_1\dots b_l)=c_1\dots c_l:=\tau\in\fS_l$.  Since $\si_1$ is a subsequence of $\pi$ there must be indices $i_1<\dots<i_l$ such that $\pi(i_j)=c_j$ for $1\le j\le l$.  Because of this and the fact that $c_1,\dots,c_l$ is a permutation of $1,\dots,l$ it must be that $\pi'_1=\tau^{-1}$.  Also $\si(j)=b_j$ for $1\le j\le l$ and $\si_1=\tau$ implies that $\tau^{-1}$ is a subsequence of $\si^{-1}$.  So $\pi'_1$ is a subsequence of $\si^{-1}$ as desired.
\eprf

There is another way to use information about one value of the antipode map to determine a second.
As in the theory of pattern avoidance, we consider the {\em diagram} of a permutation $\si\in\fS_n$ to be the set of points  
$(i,\si(i))$,  $1\le i\le n$, in the Cartesian plane.  One can then ask, 
which of the eight operations on permutations induced by the action of the dihedral group of the square preserve coefficients of the antipode?  Aside from the trivial action, there is only one other which leaves the mutiset of coefficients invariant.  Examples in $\fS_3$ and with $\si=2413$ show that the other six actions do not preserve coefficients.  Given $\si$, let $\si^o$ be the permutation whose diagram is gotten from rotating the diagram of $\si$ by $180$ degrees.  In other words if $\si=b_1\dots b_n$, then
$$
\si^o=(n+1-b_n)\dots(n+1-b_1).
$$
The next proposition now follows easily from the fact that we always have $\st(w^o)=(\st w)^o$ and so we omit the proof.
\bpr
\label{si^o}
If $\pi,\si\in\fS_n$, then
$$
[\pi]S(\si)=[\pi^o]S(\si^o).
$$
\epr

Before we start to give formulae for the antipode on specific elements of $\SSym$, we wish to recall an important connection with ribbon Schur functions.  
Consider the map $i:\NSym\ra\SSym$ defined by
\beq
\label{i(R)}
i(R_\al)=\sum_{w\in\fS_n:\ \Des( w^{-1})=D(\al)} w
\eeq
where $D$ and $\Des$ are as defined by~\ree{D(al)} and~\ree{Des}, respectively.
This is an injective Hopf algebra map~\cite[Corollary 8.14]{GR14:hac}.

We now give some explicit, cancellation-free expressions for $S(\si)$ for various specific $\si$.  We start with the identity permutation.  This can easily be derived from Theorem~\ref{S(R)}  together with $i(R_{n})=12\dots n$ and $i(R_{1^n})=n(n-1)\dots 1$.  But we prefer to give a merge-split proof.

\bpr
\label{S(eta_n)}
We have
$$
S(12\dots n)=(-1)^n(n(n-1)\dots 1).
$$
\epr
\bprf
It suffices to put a sign-reversing involution on the terms appearing in $S(12\dots n)$ whose unique fixed point is $\pi_0=n(n-1)\dots 1$ with sign $(-1)^n$.  The only shuffle set in $S(12\dots n)$ containing $\pi_0$ is the term $(-1)^n(1\shu 2\shu\dots \shu n)$ which gives us the desired fixed point.

Now take any $\pi\neq\pi_0$.  Then there must be a smallest index $i$ such that $i+1$ appears to the right of $i$ in $\pi$.  Since the numbers $1,\dots,i-1$ appear in reverse order in $\pi$, every shuffle set in $S(12\dots n)$ containing $\pi$ must be of the form 
$$
(-1)^k(1\shu 2\shu \dots \shu (i-1) \shu \si_i \shu\dots \shu \si_k)
$$ 
for some $k$.  If we are considering an appearance of $\pi$ in a shuffle set with $|\si_i|=1$, then let the involution pair it with the occurence of $\pi$ in the merged shuffle set
$$
(-1)^{k-1}(1\shu 2\shu \dots \shu (i-1) \shu\si_i \si_{i+1} \shu \si_{i+2}\shu \dots \shu \si_k).
$$
If we are considering an appearance of $\pi$ in a shuffle set with $|\si_i|>1$, then let the involution pair it with an occurence of $\pi$ in the split shuffle set
$$
(-1)^{k+1}(1\shu 2\shu \dots \shu i \shu\si_i' \shu \si_{i+1} \dots \shu \si_k),
$$
where $\si_i'$ is $\si_i$ with $i$ removed.  It is clear from the definitions that these operations are inverses and sign-reversing.
\eprf

Using the  Theorem~\ref{inverse}, the previous proposition, and the fact that $12\dots n$ is its own inverse, we can compute when $12\dots n$ appears as a term in $S(\si)$ for any $\si$.
\bco
We have
$$
[12\dots n] S(\si)=
\case{(-1)^n}{if $\si=n(n-1)\dots 1$,}{0}{else.}
$$
\eco

Using similar techniques, we can prove the following result about $\si=n\dots 21$.  It would be interesting to give general conditions under which $[\pi] S(\si)=[\pi'] S(\si')$ where the prime denotes either reflection in a vertical axis or rotation by $\pi/2$ radians.  Note that in either case $(12\dots n)'=(n\dots 21)$.
\bpr
\label{S(de_n)}
We have 
$$
S(n\dots 21) = (-1)^n (12\dots n)
$$
and 
$$
[n\dots 21] S(\si)=
\case{(-1)^n}{if $\si=12\dots n$,}{0}{else.}
$$
\epr

We are going to generalize Propositions~\ref{S(eta_n)} and~\ref{S(de_n)} to certain permutations starting with a decreasing sequence and ending with an increasing one. 
Applying the Robinson-Schensted map to these permutations outputs a pair of column superstandard tableaux of hook shape.
 First let us introduce  the notation
$$
\eta_{k,l} = k(k+1)\dots l\qmq{and} \de_{l,k}=l(l-1)\dots k
$$
with the convention that if $k>l$, then $\eta_{k,l}$ and $\de_{l,k}$ are both the empty word.  
We will further abreviate to 
$$
\eta_k=\eta_{1,k}\ \qmq{and}\de_k=\de_{n,k}
$$   
when dealing with results for $\fS_n$.
Another useful notion for applying induction is the following. If $A$ is a  set of positive integers with $|A|=n$ and $\fS_A$ is the set of permutations (linear orderings) of the elements of $A$, then we have the standardization bijection $\st_A:\fS_A\ra\fS_n$.  We then define, for any $\si\in\fS_A$,
\beq
\label{conjugate}
S(\si)=\st_A^{-1} S(\st_A(\si)),
\eeq
where $\st_A$ is extended linearly. 

We will need the following lemma which is a refinement of the well-known fact that the alternating sum of any row of Pascal's triangle (except the first) is zero.
\ble
\label{binomial}
For any $n\ge1$ we have
$$
\sum_{k=0}^n (-1)^k(\eta_k\shu \de_{k+1})=0.
$$
\ele
\bprf
It suffices to define a sign-reversing involution without fixed points on the terms appearing in the sum.  Let $v$ be a such a term and let $k$, $1\le k \le n$, be the largest integer such that $\eta_k$ is a subword of $v$.  Then $v$ appears in 
$(-1)^k(\eta_k\shu \de_{k+1})$ and, by maximality of $k$, in $(-1)^{k-1}(\eta_{k-1}\shu \de_k)$.  Since these are the only two places $v$ appears, $v$ cancels and this is true for all $v$, leaving a sum of zero.
\eprf

Our next result is the promised generalization.  Note that concatenation takes precendence over shuffle.  Note also that the following formula is cancelation free since the terms in each summand end with a different integer.
\bth
For $1\le k< n$ we have
$$
S(\de_{k,1} \eta_{k+1,n})=\sum_{j=1}^k  (-1)^{n+k+j} [\eta_{j-1}\shu(\de_{k,j+1}\shu\de_{k+2})(k+1)]j.
$$
\eth
\bprf
We will induct on $k$ where Proposition~\ref{S(eta_n)} is the base case.  We will show that the terms of 
$S(\de_{k,1} \eta_{k+1,n})$ ending in $j$, $1<j\le k$, are as given in the summation.  The cases $j=1$ and $j>k$ are similar.  Applying  the definition of an antipode and~\ree{conjugate} gives
$$
S(\de_{k,1} \eta_{k+1,n})=-\left[\sum_{i=1}^k \de_{i,1}\shu S(\de_{k,i+1}\eta_{k+1,n})
+\sum_{i=k+1}^n \de_{k,1}\eta_{k+1,i}\shu S(\eta_{i+1,n})\right].
$$
Applying induction we see that the only terms ending in $j$ will be in the first sum since $j\le k$, and that these terms must come from  the sum
$$
\sum_{i=1}^{j-1} (-1)^{n+k+j+i-1} \left(\de_{i,1}\shu[\eta_{i+1,j-1}\shu(\de_{k,j+1}\shu\de_{k+2})(k+1)]j\right).
$$
Extracting only the terms ending in $j$ from this sum and factoring out the shared expression  $v=(\de_{k,j+1}\shu\de_{k+2})(k+1)$  gives
$$
(-1)^{n+k+j} \left[\left\{\sum_{i=1}^{j-1} (-1)^{i-1} \de_{i,1}\shu\eta_{i+1,j-1}\right\}\shu v\right] j.
$$
We can now use Lemma~\ref{binomial} (reading all the words backwards for this application)  to simplify the sum to $\eta_{j-1}$.  Plugging in this as well as the value of $v$ gives that the terms ending in $j$ are exactly
$$
(-1)^{n+k+j}  [\eta_{j-1}\shu(\de_{k,j+1}\shu\de_{k+2})(k+1)]j
$$
as desired.
\eprf

Combining the previous result with  Proposition~\ref{si^o}, using the fact that 
$$
(\eta_j\shu\de_{j+1})^o = \de_{n-j,1}\shu \eta_{n-j+1,n}
$$
and reindexing gives the following corollary.
\bco
We have
$$
S(\eta_{1,k}\de_{n,k+1})=
\sum_{j=k+1}^n (-1)^{n+k+j+1} j[k(\de_{k-1,1}\shu\de_{j-1,k+1})\shu\eta_{j+1,n}].
$$
\eco

We end this section with a couple of conjectures.  To state them, we will need to extend the previous notation.
If $A$ is a set of positive integers, then we let $\eta_A$ and $\de_A$ denote the increasing and decreasing words whose elements are $A$, respectively.  Given $A\sbe[n]$ we let $\Ab=[n]-A$ be the complement of $A$ in $[n]$.  If $a\in [n]$, then we use the abbreviation $\ab=[n]-\{a\}$.  We want to consider permutations of the form $\si_A=\de_A\eta_{\Ab}$.  Note that the previous theorem deals with the case when $A=[k]$ for some $k<n$.  Now we consider what happens when $|A|=1$.  Note that when a summand in one of these expressions contains a number greater than $n$, then the expression is considered to be the empty sum.  So, for example, the first summand in the next conjecture  is empty if $n=2$.
\bcon
Let $A=\{a\}$ where  $1<a\le n$. We have
\begin{align*}
S(\si_A)&=
(-1)^{n-1} (2\shu \de_4)31+ (-1)^{n+a}((a-1)\eta_{a-2}\shu \de_{a+1})a\\[5pt]
&\hs{40pt}+\sum_{j=2}^{a-1} (-1)^{n+j} 
[((j-1)\eta_{j-2}\shu \de_{j+1})j + ((j+1)\eta_{j-1} \shu \de_{j+2})j].
\end{align*}
\econ

Here is what we believe happens when $|A|=2$ and $2\in A$.  
\bcon
Let $A=\{a,2\}$ with $2<a\le n$.  We have
\begin{align*}
S(\si_A)&=
(-1)^n[(32\shu \de_5)41 + (12\shu \de_4)3] + (-1)^{n-1}[(1\shu (3\shu \de_5)4)2]\\[5pt]
&\hs{40pt}+\sum_{j=3}^{a-1} (-1)^{n+j} 
[((j+1)21\eta_{3,j-1} \shu \de_{j+2})j - (j21\eta_{3,j-1} \shu \de_{j+2})(j+1)].
\end{align*}
\econ

Note that all of the above shuffle expressions have coefficients $\pm1$ although when they are expanded as sums of permutations, the permutations can have larger coefficients.  Can $S(\si_A)$ always be expressed in this form?

\section{The Poirier-Reutenauer Hopf algebra}

The Poirier-Reutenauer Hopf algegra~\cite{pr:aht},  $\PSym$, is a sub-Hopf algebra of the Malvenuto-Reutenauer Hopf algebra.  It has a distinguished basis indexed by standard Young tableaux.  If $\pi\in\fS$, then let $P(\pi)$ denote the insertion tableau of $\pi$ under the Robinson-Schendsted correspondence.  The basis element corresponding to a standard Young tableau $P$ is defined to be
$$
\bP=\sum_{\pi\ :\ P(\pi)=P} \pi.
$$
For example, if $P$ is the tableau given in Figure~\ref{P}, then
$$
\bP=32154+32514+32541+35214+35241.
$$
We extend this notation in the obvious way to tableaux whose entries are not necessarily $1,\dots,n$.

\begin{figure}
$$
P=\begin{ytableau}
1&4\\
2&5\\
3
\end{ytableau}
\hs{50pt}
P'=\begin{ytableau}
1&3\\
2&5\\
4
\end{ytableau}
$$
\capt{\label{P} Two standard Young tableaux, the first one superstandard, the second not.}
\end{figure}

Let $\PSym$ be the span of the $\bP$ as $P$ runs over all standard Young tableaux.  
This is a graded Hopf algebra $\PSym=\sum_{n\ge 0} \PSym_n$ where the grading is inherited from $\SSym$.  The multiplication is given by
$$
\bP\cdot\bQ=\sum_R \bR,
$$
where the sum is over all standard Young tableaux $R$ such that $P$ is a subtableau of $R$ and $Q=\st(j(R/P))$ where $j$ is  jeu de taquin and $\st$ is the standardization map applied to tableaux.  If $|P|=n$ and $Q$ is obtained by increasing all the entries of a standard  Young tableau by $n$, then it will be convenient to also define $\bP\cdot\bQ=\bP\cdot \st(\bQ)$.  For example,
\begin{align*}
\begin{ytableau}
1&2\\
3
\end{ytableau}
\cdot
\begin{ytableau}
1&2
\end{ytableau}
&=
\begin{ytableau}
1&2\\
3
\end{ytableau}
\cdot
\begin{ytableau}
4&5
\end{ytableau}
\\[10pt]
&=
\begin{ytableau}
1&2&4&5\\
3
\end{ytableau}
+
\begin{ytableau}
1&2&5\\
3&4
\end{ytableau}
+
\begin{ytableau}
1&2&5\\
3\\
4
\end{ytableau}
+
\begin{ytableau}
1&2\\
3&5\\
4
\end{ytableau}\ .
\end{align*}

To describe the coproduct, let $\pi\cong\si$ mean that $\pi$ and $\si$ are Knuth equivalent.  Then
$$
\De(\bR) = \sum_{P,Q} \st(\bP) \otimes \st(\bQ),
$$
where the sum is over all $P,Q$ whose row words satisfy $w_P w_Q \cong w_R$, or equivalently $P(w_P w_Q)=R$.  As with the product, we will sometimes not standardize $\bP$ and $\bQ$.  To illustrate, suppose 
$$
R= \begin{ytableau}
1&3\\
2
\end{ytableau}\ .
$$
Then the words in the Knuth class of $R$ are $\pi=213$ and $\pi=231$.  So to compute $\De(\bR)$ we first look at all contatenations $\pi=\pi_1\pi_2$ where $\pi_1$ and $\pi_2$ are row words of tableaux.  Putting a space between the prefixes and suffixes, we have $213 = \emp\ 213=2\ 13 = 21\ 3= 213\ \emp$ where $\emp$ is the empty word and $231 = 2\ 31 = 23\ 1$.  Note that we do not get $\emp\ 231$ since $231$ is not the row word of any tableau.  Translating from words to tableaux gives
\begin{align*}
\De(\bR) &= \emp\otimes \begin{ytableau} 1&3\\ 2 \end{ytableau}\ +\ 
\begin{ytableau} 2 \end{ytableau}\otimes \begin{ytableau} 1&3 \end{ytableau}\ +\
\begin{ytableau} 1\\ 2 \end{ytableau} \otimes \begin{ytableau} 3   \end{ytableau}\ +\
\begin{ytableau} 1&3\\ 2 \end{ytableau} \times\emp
\\[10pt]
&\qquad  +\ \begin{ytableau} 2 \end{ytableau}\otimes \begin{ytableau} 1\\ 3 \end{ytableau}\ +\
\begin{ytableau} 2&3 \end{ytableau} \otimes \begin{ytableau} 1   \end{ytableau}
\\[10pt]
 &= \emp\otimes \begin{ytableau} 1&3\\ 2 \end{ytableau}\ +\ 
\begin{ytableau} 1 \end{ytableau}\otimes \begin{ytableau} 1&2 \end{ytableau}\ +\
\begin{ytableau} 1\\ 2 \end{ytableau} \otimes \begin{ytableau} 1   \end{ytableau}\ +\
\begin{ytableau} 1&3\\ 2 \end{ytableau} \times\emp
\\[10pt]
&\qquad  +\ \begin{ytableau} 1 \end{ytableau}\otimes \begin{ytableau} 1\\ 2 \end{ytableau}\ +\
\begin{ytableau} 1&2 \end{ytableau} \otimes \begin{ytableau} 1   \end{ytableau}\ .
\end{align*}

The antipode in $\PSym$ seems to be even more complicated than the one in $\SSym$.  But we at least have a conjecture for certain hook-shaped tableaux.  Let $\la$ be a partition.  The {\em column superstandard Young tableau of shape $\la$}, $P_\la$, is obtained by filling the first column with the numbers $1,\dots,k$, then the second column with the numbers $k+1,\dots,k+l$, and so forth.  In Figure~\ref{P},  we have $P=P_{(2,2,1)}$ while $P'$ is not column superstandard.  Recall that the {\em descent set} of a standard Young tableau $P$, $\Des P$, is the set of all $i$ such $i+1$ is in a lower row.  From properties of the Robinson-Schensted correspondence it follows that $P(\pi)=P$ implies $\Des\pi^{-1}=\Des P$.  If $\la=(n,1^k)$ is a hook and $P=P_\la$  then the converse is also easily  seen to be true.  It follows that $\bP_\la=i(R_{1^k,n})$ where $i$ is the map in~\ree{i(R)}.  Together with Theorem~\ref{S(R)}, this proves the following result.

\bth
\label{PSymThm}
If $\la$ is a hook, then in $\PSym$ we have
$$
S(\bP_\la)=(-1)^{|\la|} \bP_{\la^t}.
$$
\eth

We have not been able to give an involution proof of this result except in some special cases.  We have already done this  when $\la$ is a single row; see Proposition~\ref{S(eta_n)}.   It is also possible to use this technique on two-row hooks.
\bpr
If $\la=(n-1,1)$, then in $\PSym$ we have
$$
S(\bP_\la)=(-1)^{|\la|} \bP_{\la^t}.
$$
\epr
\bprf
First note that the Knuth class for $P_\la$ consists of the permutations 
$$
\pi_i = \eta_{2,i+1} 1 \eta_{i+2,n},
$$
where $1\le i\le n-1$.  So the terms in Takeuchi's expansion for $S(\bP_\la)$ are of the form $(-1)^k \bP_1\cdot\ \dots\ \cdot \bP_k$ where the $\bP_j$ come from a concatenation of subwords of some $\pi_i$ which are row words of tableaux.  The variable $k$ will always denote the number of factors.  We will associate with such a term the pair $(i,\al)$ where $\pi_i$ is the permutation giving rise to the term and  $\al=(\al_1,\dots,\al_k)$ is a composition where $\al_j=|P_j|$ for all $j$.  Note that this pair fully determines the corresponding term.   Finally, we will not standardize the $P_j$ but rather write $\st(P_j)$ explicitly if we need to do so.

Initially, we will be cancelling  a single term in the expansion with a pair of terms which have the same sum but with opposite sign.  The pair will be obtained by two different mergings of the single term.  To describe the involution, it will be simplest to describe a bijection $f:D\ra R$ whose domain elements are certain single terms and whose range elements  are pairs of terms.  Merging will correspond to applying $f$ while splitting will correspond to applying $f^{-1}$.

Our first function $f_1:D_1\ra R_1$ has a domain all terms with 
$$
\bP_j = \raisebox{-5pt}{\begin{ytableau} 1 \end{ytableau}}
$$
for some $j\neq k$.  Note that it is also not possible to have $j=1$ since none of the $\pi_i$ begin with a one.  Thus the associated pairs for these terms are of the form $(i,\al)$ where $i\neq n-1$ and $\al_j=1$.  The summands in the range will consist of all terms where $1$ appears in a $\bP_j'$  such that  $j\neq 1$ and there is at least one other element in the tableau $P_j'$.  
Write $\bP_1\cdot\ \dots\ \cdot \bP_k = A\cdot \bP_j  \cdot \bP_{j+1}\cdot B$ and let
$$
f_1(A\cdot \bP_j  \cdot \bP_{j+1}\cdot B) = A\cdot \bP_j' \cdot B + A\cdot \bP_j''\cdot B,
$$
where the first summand is determined by the pair $(i,\al')$ and the second by $(i+1,\al')$ where $\al'$ is $\al$ with the parts $\al_j=1$ and $\al_{j+1}$ replaced by a single part $\al_j'=1+\al_{j+1}$.  For example, if $n=5$, $i=2$, and $\al=(1,1,1,2)$, then $\pi_2=23145$ and the term corresponding to $\al$ is 
$$
\bP_1\cdot\bP_2\cdot\bP_3\cdot\bP_4
=
 \raisebox{-5pt}{\begin{ytableau} 2 \end{ytableau}}
\cdot
 \raisebox{-5pt}{\begin{ytableau} 3 \end{ytableau}}
\cdot
 \raisebox{-5pt}{\begin{ytableau} 1 \end{ytableau}}
\cdot
 \raisebox{-5pt}{\begin{ytableau} 45 \end{ytableau}}\ .
$$
Now $j=3$ and so $\al'=(1,1,1+2)=(1,1,3)$.  So $f$ maps this product to the sum of the products associated with $\al'$ in $\pi_2$ and $\pi_3=23415$ which is 
$$
f_1(\bP_1\cdot\bP_2\cdot\bP_3\cdot\bP_4)
=
 \raisebox{-5pt}{\begin{ytableau} 2 \end{ytableau}}
\cdot
 \raisebox{-5pt}{\begin{ytableau} 3 \end{ytableau}}
\cdot
 \raisebox{-5pt}{\begin{ytableau} 1 & 4 & 5 \end{ytableau}}
+
\raisebox{-5pt}{\begin{ytableau} 2 \end{ytableau}}
\cdot
 \raisebox{-5pt}{\begin{ytableau} 3 \end{ytableau}}
\cdot
 \raisebox{-5pt}{\begin{ytableau} 1 & 5\\ 4 \end{ytableau}}\ .
$$

We first need to verify that $f$ is well defined in that the subwords of $\pi_i$  and $\pi_{i+1}$ defined by $\al'$ are indeed row words of tableaux.  Since $1$ must be the subword of $\pi_i$ corresponding to $\al_j=1$, the word corresponding to $\al_{j+1}$ must be $\eta_{i+2,l}$ for some $l$.  It follows that the subwords of $\pi_i$ and $\pi_{i+1}$ corresponding to $\al_j'=1+\al_{j+1}$ are $1\eta_{i+2,l}$ and $(i+2)1\eta_{i+3,l}$, respectively.  It is easy to see that both of these are row words.

Next we need to show that a product and its image under $f_1$ cancel each other out in Takeuchi's expansion.  Clearly the initial product and the two image products are of opposite sign.  So it suffices to show that $\bP_j\cdot\bP_{j+1}=\bP_j'+\bP_j''$.  Using the description of the corresponding subwords in the previous paragraph gives
\begin{align*}
\bP_j\cdot\bP_{j+1}
&= 
 \raisebox{-5pt}{\begin{ytableau} 1 \end{ytableau}}
\cdot
\barr{|c|c|c|c|}\hline \rule{0pt}{14pt} i+2 & i+3 & \dots & l \\ \hline\earr
\\[10pt]
&=
\barr{|c|c|c|c|c|}\hline \rule{0pt}{14pt} 1& i+2 & i+3 & \dots & l \\ \hline\earr\
+
\raisebox{-10pt}{
$\barr{|c|c|c|c|}\hline \rule{0pt}{14pt} 1 & i+3 & \dots & l 
\\ \hline
\rule{0pt}{14pt} i+2
\\ \cline{1-1}
\earr$
}
\\[10pt]
&=
\bP_j'+\bP_j''.
\end{align*}

Finally, we need to prove that $f_1$ is a bijection.  We do this by constructing its inverse.  Suppose we are given a product such that the factor $\bP_j'$ containing one satisfies $j\neq 1$ and  $|P_j'|\ge2$.  We must find the product to add to the given one so that it can be mapped back to the domain.  Suppose the product is associated with a pair $(i,\al')$.  Then considering the subwords of $\pi_i$ which contain $1$ and which are row words of tableaux, we see that there are only two possibilities for $P_j'$, namely
$$
P_j'=\barr{|c|c|c|c|c|}\hline \rule{0pt}{14pt} 1& i+2 & i+3 & \dots & l \\ \hline\earr\
\qmq{or}
\raisebox{-10pt}{
$\barr{|c|c|c|c|}\hline \rule{0pt}{14pt} 1 & i+2 & \dots & l 
\\ \hline
\rule{0pt}{14pt} i+1
\\ \cline{1-1}
\earr$
}
$$
for some $l$. If  $P_j'$ is the first (respectively, second) of these tableaux,  then we pair the given product with the product associated with the pair $(i+1,\al')$ (respectively, $(i-1,\al')$).  It is now easy to verify that adding this pair simplifies to a single product which was mapped to the pair by $f_1$.

The next part of the involution is similar to the first, so we will only provide a description of the map.  Define a map 
$f_2:D_2\ra R_2$ where the domain contains all products with
$$
\bP_1= \raisebox{-5pt}{\begin{ytableau} 2 \end{ytableau}}
\qmq{and}
\bP_k = \raisebox{-5pt}{\begin{ytableau} 1 \end{ytableau}}\ .
$$
Note that this forces the term to be associated with the pair $(n-1,\al)$ where $\al_1=\al_k=1$.  Let $f_2$ map this to the sum of the terms associated with the pairs $(1,\al')$ and $(n-1,\al')$ where $\al'$ is $\al$ with its first two parts replaced by $\al_1'=1+\al_2$.  One can then verify that the range contains sums of all pairs where the first (respectively second) summand contains one and two (respectively, two and three) in the first product tableau and $n$ (respectively, $1$) as a singleton in the last.

From the descriptions of $D_1,D_2,R_1,R_2$ and an examination of which subwords of the $\pi_i$ can be row words, we see that the only terms which remain uncancelled so far are those associated with pairs of the form $(1,\al)$ where $\al_1,\al_k\ge2$.  Among these products, the only one which can produce $\bP_{\la^t}$ is the one associated with $\al=(2,1,1,\dots,1,2)$.  So we will be finished if we can define a sign-reversing involution on the tableaux $P\neq P_{\la^t}$ which occur in expanding the products under consideration.

Since all our products come from $\pi_1$, it suffices to specify $\al$ to focus on a given product whose expansion contains $\bP$.  So suppose the product corresponding to $\al$ has $\bP$ as a term.   Let $l$ be the largest integer such that $1,\dots,l$ are in the same cells in both $P$ and $P_{\la^t}$.  Equivalently we have $1,\dots,l$ in the first column of $P$ but $l+1$ is not.
Since $P$ comes from a product associated with $\pi_1$ and $\al_1\ge2$ we have $l\ge2$.  And since $P\neq P_{\la^t}$, $l\le n-2$.   Let $\bP_j$ be the  factor in the product containing $l$.  If $l+1\not\in \bP_j$, then $\bP$ can be cancelled with its appearance in the product corresponding to $\al'$ obtained from $\al$ by replacing $\al_j$ with $\al_j+\al_{j+1}$.  If $l+1\in\bP_j$, then form $\al'$ by splitting $\al_j=\al_j'+\al_{j+1}'$ where $\al_j'$ corresponds to the prefix (respectively, suffix) of the $\al_j$ subword of $\pi_1$ consisting of those elements less than or equal to (respectively, greater than) $l$.  Note that the bounds on $l$ guarantee that $\al'$ will still satisfy $\al_1',\al_l'\ge 2$.  To illustrate, suppose $n=5$ and $\al=(2,3)$.  Then $\pi_1=21345$ and so we are considering the product
$$
\begin{ytableau}
1 \\ 2 
\end{ytableau}\
\cdot\
\begin{ytableau}
3&4&5
\end{ytableau}\
=\
\begin{ytableau}
1 & 3& 4 & 5\\ 2 
\end{ytableau}\
+\
\begin{ytableau}
1 & 4 & 5\\ 2 \\ 3
\end{ytableau}\ .
$$
Comparing the first summand to 
$$
\bP_{\la^t} =\
\begin{ytableau}
1 & 5 \\ 2 \\ 3 \\ 4
\end{ytableau}
$$
we see that $l=2$.  Since $l+1=3$ does not occur in the same tableau as $2$ in the product, we merge and cancel this term with the tableau corresponding to $\al'=(5)$ which is
$$
-\
\begin{ytableau}
1 & 3& 4 & 5\\ 2 
\end{ytableau}\ .
$$
For the second summand we have $l=3$ which is in the same tableau as $l+1=4$ in the product.  So we split $\al$ into $\al'=(2,1,2)$ and see that this $\bP$ will cancel with one of the terms in
$$
-\
\begin{ytableau}
1\\ 2 
\end{ytableau}\
\raisebox{5pt}{$\cdot$}\
\begin{ytableau}
3
\end{ytableau}\
\raisebox{5pt}{$\cdot$}\
\begin{ytableau}
4 & 5 
\end{ytableau}\ .
$$

The proof that this is a well-defined, sign-reversing involution is similar to previous arguments we have seen earlier and so is omitted.  This completes the demonstration of the proposition.
\eprf

\section{Future work and open problems}
\label{fwo}

We hope that this article will just be a first step in the exploration of the use of sign-reversing involutions to derive formulas for antipodes.  In addition to the conjectures and questions already raised, here are three directions which would be interesting to explore.

\ben
\item  Are there other Hopf algebras where the split or merge idea can be used to derive nice, preferably cancellation-free, formulas for $S$?  Even more ambitious, is there a (meta)-involution which can be used to prove antipode identities for many different Hopf algebras at once?  As was noted at the end of Section~\ref{hag}, one can use an involution for $\bbF[x]$ which is the special case of the one for $\cG$ where the graph has no edges.  We should also mention that the proof in Section~\ref{hag} is based on discussions with  Nantel Bergeron.  The involution we introduce is closely related to the one of Bergeron and  Ceballos~\cite{BC15:has} for their Hopf algebra of subword complexes.  A formula for the antipode of a Hopf algebra of abstract simplicial complexes has been computed using these techniques by the first author in joint work with Hallam and Machacek~\cite{bhm:cha}.  Samantha Dahlberg (private communication) has used our method to compute the antipode for a Hopf algebra on involutions.  Finally, Eric Bucher and Jacob Matherne (private communication) have used the merge-split technique to determine the antipode for the restriction-contraction Hopf algebra on uniform matroids.

\item  Can one obtain a full cancellation-free formula for the antipode in $\NSym$ in the immaculate basis?  We attempted to at least do the three-row case, but the expressions in terms of frozen tableaux became increasingly complicated.  However, there may be some other idea which is needed to unify all the cases.
\item  We know from equation~\ree{S(s)} that the antipode for $\Sym$ is  particularly simple when expressed in terms of the Schur basis.  Is there a way to derive this beautiful formula using a sign-reversing involution?
\een

{\em Acknowledgments.}
We are grateful to Nantel Bergeron and Darij Grinberg for enlightening discussions.  Grinberg also found a proof of Theorem~\ref{PSymThm} while it was still a conjecture in a previous version of this paper.  Finally,, we thank the two anonymous referees for suggestions which improved the manuscript.  In particular, one referee pointed out the connection of some of our results with ribbon Schur funcitons which was very useful.

\bibliographystyle{alpha}

\newcommand{\etalchar}[1]{$^{#1}$}

\end{document}